\documentclass[preprint,10.5pt]{elsarticle}
\pdfoutput=1

\usepackage[margin=1in]{geometry}

\usepackage{tikz}
\usepackage{array}
\usepackage{amsmath}
\usepackage{amssymb}
\usepackage{booktabs}
\usepackage{subcaption}
\usepackage{soul}
\usepackage{xcolor}
\usepackage{subfig}
\usepackage{multirow}
\usepackage[export]{adjustbox}
\usepackage{ulem}
\usepackage{algorithm}
\usepackage{algpseudocode}
\usepackage{todonotes}
\usepackage{nicematrix}

\begin{document}

\begin{frontmatter}

\title{EKF--SINDy: Empowering the extended Kalman filter with sparse identification of nonlinear dynamics}

\author[1]{Luca Rosafalco\corref{cor1}}
\ead{luca.rosafalco@polimi.it}
\author[1]{Paolo Conti}
\ead{paolo.conti@polimi.it}
\author[2]{Andrea Manzoni}
\ead{andrea1.manzoni@polimi.it}
\author[1]{Stefano Mariani}
\ead{stefano.mariani@polimi.it}
\author[1]{Attilio Frangi}
\ead{attilio.frangi@polimi.it}
\address[1]{Department of Civil and Environmental Engineering,
Politecnico di Milano \\ Piazza L. da Vinci 32,
20133 - Milano (Italy)}
\address[2]{MOX, Department of Mathematics,
Politecnico di Milano \\ Piazza L. da Vinci 32,
20133 - Milano (Italy)}

\cortext[cor1]{Corresponding author}


\begin{abstract}
Measured data from a dynamical system can be assimilated into a predictive model by means of Kalman filters. Nonlinear extensions of the Kalman filter, such as the Extended Kalman Filter (EKF), are required to enable the joint estimation of (possibly nonlinear) system dynamics and of input parameters. To construct the evolution model used in the prediction phase of the EKF, we propose to rely on the Sparse Identification of Nonlinear Dynamics (SINDy). SINDy enables to identify the evolution model directly from preliminary acquired data, thus avoiding possible bias due to wrong assumptions and incorrect modelling of the system dynamics. Moreover, the numerical integration of a SINDy model leads to great computational savings compared to alternate strategies based on, e.g., finite elements. Last, SINDy allows an immediate definition of the Jacobian matrices required by the EKF to identify system dynamics and properties, a derivation that is usually extremely involved with physical models. As a result, combining the EKF with SINDy provides a data-driven computationally efficient, easy-to-apply approach for the identification of nonlinear systems, capable of robust operation even outside the range of training of SINDy. To demonstrate the potential of the approach, we address the identification of a linear non-autonomous system consisting of a shear building model excited by real seismograms, and the identification of a partially observed nonlinear system. The challenge arising from the use of SINDy when the system state is not entirely accessible has been relieved by means of time-delay embedding. The great accuracy and the small uncertainty associated with the state identification, where the state has been augmented to include system properties, underscores the great potential of the proposed strategy, paving the way for the setting of predictive digital twins in different fields.
\end{abstract}


\end{frontmatter}

\section{Introduction}
\label{sec:introduction}

Dynamical system identification is crucial to enable the construction of predictive digital twins \cite{art:Wagg20}, and to implement control strategies for engineering systems with applications to condition-based maintenance of civil structures \cite{art:Matt24}, or to fuel consumption reduction in air transportation \cite{art:Guastoni23}. Despite advances in Machine Learning (ML) and Data Science \cite{art:Jordan15}, identifying explainable, reduced dimensional models of dynamic processes from big data is an open field of research \cite{art:Brunton16}. Searching for the physical relations that underlie a certain dynamic process is probably the best way to obtain robust and generalisable models, a characteristic uncommon in most of ML techniques. Consequently, making reliable predictions in scenarios lacking collected data is generally challenging or impossible \cite{art:Olivier21}. This type of uncertainty is commonly referred to as epistemic uncertainty. To address this challenge, we employ the Sparse Identification of Nonlinear Dynamics (SINDy) proposed in \cite{art:Brunton16}. SINDy constructs robust and generalisable models by assuming that only a few important terms govern the dynamics of the considered system. This assumption holds for many physical systems when an appropriate basis is used to describe the space of functions governing their dynamics \cite{art:Brunton16}.

\begin{figure}[b!]
\centering
\includegraphics[width=\linewidth]{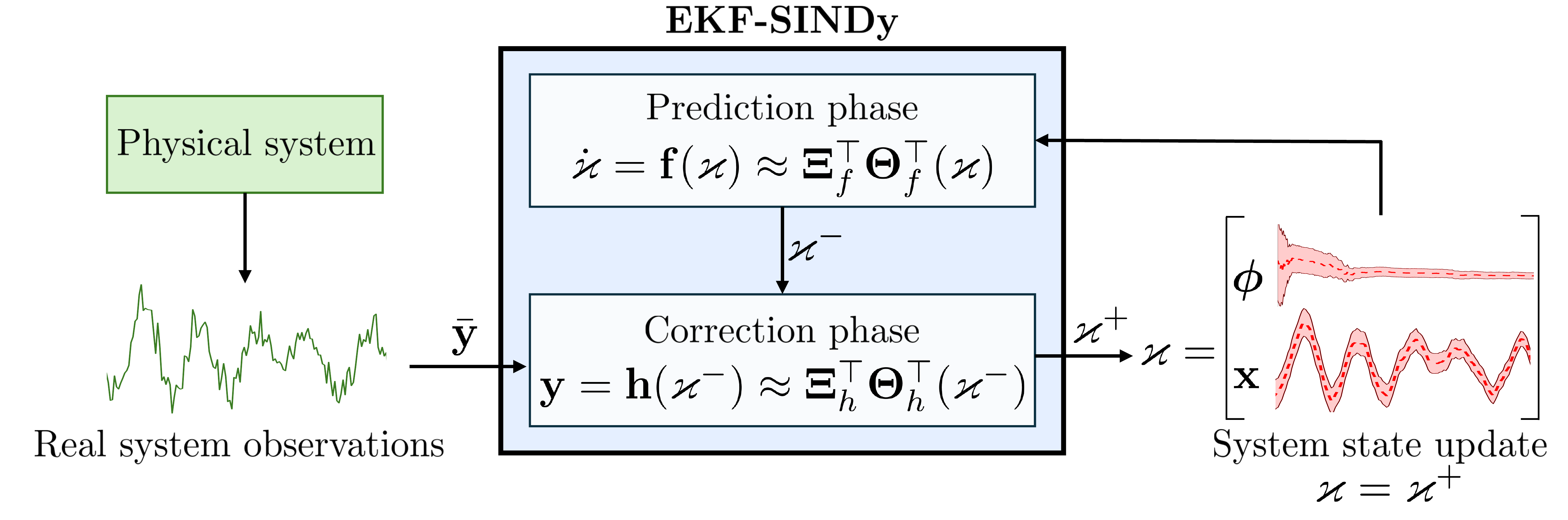} 
\caption{\footnotesize EKF-SINDy methodology. In the prediction phase, the dynamical system identified by SINDy, together with random walk equation modelling the parameter evolution, advances the augmented state $\boldsymbol{\varkappa}=[\mathbf{x}^{\top},\boldsymbol{\phi}^{\top}]^{\top}$—comprising system state $\mathbf{x}$ and parameters $\boldsymbol{\phi}$—to a prior estimate $\boldsymbol{\varkappa}^{-}$; the correction phase involves the EKF assimilating system observations $\bar{\mathbf{y}}$, resulting in the updated augmented state $\boldsymbol{\varkappa}^{+}$.}
\label{fig:graph_abs_ekf_sindy}
\end{figure}

To construct a digital twin, the model must evolve over time to reflect potential changes in its physical counterpart. This evolution should be, possibly online, driven by real-world data \cite{art:Wagg24}. In this work, we have exploited an Extended Kalman filters (EKF) to perform data assimilation \cite{art:Kalman60} leveraging a prediction-correction scheme. During the prediction stage, the SINDy model evolves the system state. The procedure can also accommodate potential external forcing rendering the system not autonomous. During the correction stage, acquired data are used to refine previous estimates, update the physical parameters governing the system dynamics, and potentially reduce the associated uncertainties. As a result, a joint estimation of the system is achieved (termed \textit{joint} because both the state of the system and the physical parameters underlying the dynamics are jointly updated, see \cite{proc:Wan00}). A schematic representation of the EKF-SINDy procedure is reported in Fig. \ref{fig:graph_abs_ekf_sindy}. More details on this methodology will be discussed in Sec. \ref{sec:methodology}.

In realistic applications, the observation variables assimilated by the EKF do not necessarily match up with the system state on which SINDy is performed. This mismatch arises in two scenarios: where the observations are fewer than the system state variables, potentially limiting the full description of the system dynamics, and where observations are in excess, potentially overwhelming the EKF-SINDy framework with redundant information. In the former case, we can resort to uplifting techniques such as, e.g., embedding techniques \cite{art:Brunton17, bakarji2023discovering} or recurrent decoder networks \cite{williams2023sensing}, to recover the hidden state components from the available, observed time-series. To address the latter case, dimensionality reduction techniques, employing proper orthogonal decomposition and/or encoder neural networks, have been coupled with SINDy to simultaneously reduce dimensionality while learning the corresponding governing dynamical system as in \cite{champion2019data, bakarji2023discovering, conti2023reduced}. In this work, we deal with the first scenario, focusing on revealing hidden variables when only partial observations are available. The integration in the method of dimensionality reduction techniques to deal with high-dimensional data will be the focus of future investigations.

In the recent literature, other works employ a ML-based identified model in the prediction phase of the KF. In \cite{proc:Coskun}, predictions were made using recurrent Neural Networks (NN), specifically Long Short-Term Memory (LSTM) modules; LSTM modules were also utilised to model the process and observation noise required by the filter. The linear KF was used to estimate only the dynamics of the observed system. In contrast, our work explicitly calibrates the physical parameters affecting the dynamical process, like in \cite{art:Liu24} where NNs and EKF were combined to perform joint estimation of mechanical systems. This work improved upon previously proposed variational autoencoders-based approaches, such as the deep Markov models \cite{proc:Krishnan17}, by relying on the EKF equation to infer the latent variables governing the system dynamics. Compared to \cite{art:Liu24}, our methodology looks better suited for a physical process whose dynamics is controlled by a few important terms. This is because NNs require large volume of training data, do not easily account for known symmetries and constraints, and struggle to handle epistemic uncertainties \cite{art:Kaiser18}.

SINDy has previously been coupled with the KF in \cite{proc:Wang22}. In that work, the estimation procedure directly targeted the SINDy coefficients. However, as previously mentioned, relying on the linear KF precludes the possibility of performing joint estimation. At variance with \cite{proc:Coskun}, the proposal of \cite{proc:Wang22} focused on identifying the coefficients of dynamic models without correcting the state predictions.

The application of the proposed approach to experimental data appears to be extremely promising. Indeed, SINDy can be directly applied to preliminary acquired data, thereby mitigating the risk of introducing a bias when modelling the system dynamics. If data driven model identification techniques like SINDy are not used, this bias can negatively impact the filtering outcome. In such cases, other possible mitigation strategies, for example based on the recently proposed integration with echo state networks, must be considered \cite{art:Novoa24}.

In this study, we have tested our methodology using simulated data for two primary reasons: first, to showcase the method potential and its ability to handle non autonomous and partially observed nonlinear systems; moreover, to demonstrate that the proposed methodology may be advantageous even when a physics-based model of the system is available. In the latter case, the model identified by SINDy serves as a surrogate for the physics-based model enabling significant computational time savings compared to the integration of nonlinear models. These computational savings are relevant to the online application of the procedure, in which data assimilation is performed by the EKF using a trained SINDy model to evolve the system dynamics. The training of the SINDy model, on the other hand, is carried out in an offline phase, which may be computationally intensive.

Relying on fast computational approaches has recently proven crucial in many applications, such as the design optimisation of Micro-Electro Mechanical Systems (MEMS) \cite{art:Vizzaccaro24}, and inverse analysis in structural health monitoring \cite{art:CAS21,chp:CSAI22}, just to mention two relevant cases dealing with structural dynamics, both at the micro- and the macro-scale. The need of computational efficiency also underlies the use of SINDy in control applications \cite{art:Kaiser18,art:abdullah23}. Alternatively, a recent work proposed to reduce the computational burden of the estimation process by using a rank-reduced version of the KF \cite{proc:Schmidt23}. Another major advantage of using SINDy is the greatly simplified computation of the Jacobian matrices required by the EKF formulation, compared to the involved derivation typically encountered even for relatively simple mechanical systems \cite{art:RUENG24}.

The remainder of the paper is organized as follows. The methodology is presented in Sec. \ref{sec:methodology}: first, KF is introduced; then the rationale behind SINDy is discussed; finally, the application of EKF-SINDy to joint estimation of autonomous and non-autonomous dynamical systems is detailed. Two case studies are presented in Sec. \ref{sec:results}. In Sec. \ref{sec:shearBuildingResults}, a first case discusses the identification of a shear building subjected to a seismic event. Real seismograms are used as excitations, demonstrating the capability of our procedure of handling non-autonomous systems. A second case study, presented in Sec. \ref{sec:resonatorResults}, focuses on the identification of a partially observed nonlinear resonator. To address the impossibility to observe the whole system state, we precede the training of the SINDy model by applying the time-delay embedding \cite{art:Champion19}. Final considerations are collected in Sec. \ref{sec:conclusions}, along with a discussion of future developments.
The source code of the proposed method is made available in the public repository \texttt{EKF-SINDy} \cite{EKFSINDY_repo}.

\section{Methodology}
\label{sec:methodology}

\subsection{Extended Kalman Filtering}

Considering a dynamical system of interest, Kalman filters exploit a state--space representation of the type:
\begin{equation}
    \dot{\mathbf{x}}(t) = \frac{\text{d}}{\text{d} t}\mathbf{x}(t) = \mathbf{f}(\mathbf{x}(t)),
    \label{eq:dynamicSystem}
\end{equation}
where: $\mathbf{x}(t)\in\mathbb{R}^n$ is the state vector of the system at time $t$; $\mathbf{f}:\mathbb{R}^n\rightarrow \mathbb{R}^n$ is the function of $\mathbf{x}(t)$ describing the dynamics of the system. Nonlinear versions of the KF, like the EKF, are required if $\mathbf{f}$ is nonlinear. Hereon, we rely on the extended Kalman filter (EKF); the reader may refer, e.g., to \cite{book:Simon06_13} for a complete derivation and analysis of the EKF theory.

Dynamical systems can be used to model the response of a building to seismic excitations, or the behavior of MEMS \cite{book:Corigliano4}. In general, $\mathbf{f}$ is not exactly known. Thus, it becomes important to assimilate incoming data of the system to update the model and the resulting predictions. For instance, these data can consist of floor acceleration measurements for a building \cite{art:Farrar98}, or displacement measurements for capacitive sensors. According to the discrete acquisition of system measurements, 
data assimilation is performed at discrete $t_j$ instants of observation,
with $j=1,\ldots,T$. As state variables may not be directly compared against incoming data, the following observation equation:
\begin{equation}
\mathbf{y}(t_j)=\mathbf{h}(\mathbf{x}(t_j))
    \label{eq:obsEqDiscrete}
\end{equation}
is usually introduced to encode the state-to-observation map. Specifically: $\mathbf{h}:\mathbb{R}^n\rightarrow \mathbb{R}^o$ is a possibly nonlinear function; $\mathbf{y}(t_j)\in\mathbb{R}^{o}$ are the quantities to be compared with the incoming observations $\bar{\mathbf{y}}(t_j)\in\mathbb{R}^{o}$. According to Eqs. \eqref{eq:dynamicSystem} and \eqref{eq:obsEqDiscrete}, model predictions are performed employing a continuous time description, while data assimilation is carried out at discrete time steps. For this reason, a hybrid (continuous--discrete) formulation of the KF will be considered in the following.

A prediction--correction scheme is then adopted to perform online data assimilation. In the prediction phase, the state of the system is advanced in time through $\mathbf{f}$. This prediction is then corrected through the comparison with incoming data. 
Hereon, we will use the superscript ``$-$" to refer to quantities not yet updated by the correction stage, and the superscript ``$+$" for quantities updated by the correction stage.

To perform the prediction--correction scheme, state variables are treated as random variables. Therefore, Eqs. \eqref{eq:dynamicSystem} and \eqref{eq:obsEqDiscrete} are modified as follows:
\begin{subequations}
    \begin{align}
    \dot{\mathbf{x}}(t) &= \mathbf{f}(\mathbf{x}(t)) + \mathbf{q}_x,
    \label{eq:dynamicSystem_noise}\\
    \mathbf{y}(t) &= \mathbf{h}(\mathbf{x}(t)) + \mathbf{r},
    \label{eq:obs_noise}
\end{align}
\end{subequations}
where: $\mathbf{q}_x$ is the process noise vector, assumed to be sampled from a white, zero mean stochastic process $\mathbf{q}_x\sim \mathcal{N}(\mathbf{0},\mathbf{Q}_x)$ with diagonal covariance matrix $\mathbf{Q}_x\in\mathbb{R}^{n\times n}$; $\mathbf{r}\in\mathbb{R}^o$ is the observation noise vector, assumed to be sampled from a white, zero-mean stochastic process $\mathbf{r}\sim \mathcal{N}(\mathbf{0},\mathbf{R})$ with diagonal covariance matrix $\mathbf{R}\in\mathbb{R}^{o\times o}$. The vector $\mathbf{q}_x$ accounts for uncertainties in the mapping $\mathbf{f}$, while $\mathbf{r}$ accounts for the noise affecting the incoming data. Their trade-off determines how much the filter relies on model predictions with respect to the acquired data \cite{book:Simon06_5}.

The Kalman filter follows the system evolution by looking at $\hat{\mathbf{x}}$ and $\mathbf{P}$, where: $\hat{\mathbf{x}}=\mathbb{E}[\mathbf{x}]$ is the mean value of $\mathbf{x}$; $\mathbb{E}$ computes the expected value of a quantity; $\mathbf{P}=\mathbb{E}[(\mathbf{x}-\hat{\mathbf{x}})(\mathbf{x}-\hat{\mathbf{x}})^\top]$ is the covariance matrix associated to $\mathbf{x}$. This is equivalent to assume that $\mathbf{x}$ features a Gaussian distribution. This assumption is quite natural in system identification. A Gaussian probability distribution is associated with $\mathbf{x}$ and $\mathbf{y}$ because it is the distribution with maximum information entropy given the specified mean and variance \cite{book:Yuen10_2}.

\paragraph{Prediction stage}

The values of $\hat{\mathbf{x}}^-(t_j)$ and $\mathbf{P}^-(t_j)$ at $t_j$ are predicted starting from $\hat{\mathbf{x}}^+(t_{j\!-\!1})$ and $\mathbf{P}^+(t_{j\!-\!1})$ at $t_{j\!-\!1}$. These values have been obtained by correcting the previous estimates with data acquired at $t_{j-1}$. The prediction at $t_j$ is performed as in the following:
\begin{subequations}
    \begin{align}        \hat{\mathbf{x}}^-(t_j) = \hat{\mathbf{x}}^+(t_{j\!-\!1})+\int_{t_{j-1}}^{t_j} \mathbf{f}(\hat{\mathbf{x}}^-(\tau))  \text{d}\tau, \label{eq:meanTransition}\\
    \mathbf{P}^-(t_j) = \mathbf{P}^+(t_{j-1})+\int_{t_{j\!-\!1}}^{t_j} \mathbf{F}^-(\tau)\mathbf{P}^-(\tau) + \mathbf{P}^-(\tau) \mathbf{F}^-(\tau)^\top + \mathbf{Q}_x \text{ }\text{d}\tau,
    \label{eq:covTransition}
    \end{align}
    \label{eq:Transition}
\end{subequations}
where $\mathbf{F}\in\mathbb{R}^{n\times n}$ is the Jacobian of $\mathbf{f}(\mathbf{x}(t))$. The matrix $\mathbf{F}^-(\tau)$ is evaluated as:
\begin{equation}
 \mathbf{F}^-(\tau) =\frac{\partial \mathbf{f}}{\partial \mathbf{x}^{\top}}\bigg|_{\hat{\mathbf{x}}^-(\tau)}.
 \label{eq:Jacobian_def}
\end{equation}

Different methods can be used to integrate Eq. \eqref{eq:Transition}. For its computational efficiency, we have considered the Euler forward method leading to:
\begin{subequations}
    \begin{align}        
    \hat{\mathbf{x}}^-(t_j) = \hat{\mathbf{x}}^+(t_{j\!-\!1})+\Delta t ~\mathbf{f}(\hat{\mathbf{x}}^+(t_{j\!-\!1})), \label{eq:meanTransition_Euler}\\
    \mathbf{P}^-(t_j) = \mathbf{P}^+(t_{j\!-\!1})+\Delta t \left(\mathbf{F}^+(t_{j\!-\!1})\mathbf{P}^+(t_{j\!-\!1}) + \mathbf{P}^+(t_{j\!-\!1}) \mathbf{F}^{+^\top}(t_{j\!-\!1}) + \mathbf{Q}_x\right),
    \label{eq:covTransition_Euler}
    \end{align}
\end{subequations}
where $\mathbf{F}^{+}(t_{j-1})$ is evaluated for $\hat{\mathbf{x}}^+_{j-1}$ at $t_{j-1}$. In general, $\mathbf{Q}_x$ may modify in time; however, we neglect such possibility in this work. Moreover, we assume a constant time step $\Delta t = t_j - t_{j\!-\!1}$.

The Euler forward method is a first-order explicit time integration technique, meaning that the local truncation error generated scales linearly with the time step $ \Delta t$, provided this latter is chosen to be sufficiently small. This time integration method has been selected for its simplicity. However, more accurate methods can be employed, such as, e.g., higher-order Runge--Kutta methods, as the selection of the time integration technique does not affect the generality of the proposed procedure. In the next section we will show that the local truncation error performed by the Euler forward method does not preclude good outcomes for the system estimation, also thanks to the adopted prediction--correction strategy. Similarly, implicit time stepping schemes can also be considered in order not to strongly constrain the choice of the time step, but it implies the solution of a (non)linear system of equations at each time step, in the case of (non)linear dynamical systems. If the employed integration algorithm is not unconditionally stable, as typical of explicit time integration schemes, preliminary tests must be performed to verify that $\Delta t$ does not induce any loss of convergence in the integration of system dynamics.

\paragraph{Correction stage}

Once the prediction stage is accomplished, the correction stage is performed by assimilating the incoming data. Specifically, KF performs corrections in order to minimise the trace of $\mathbf{P}$. In the EKF, the following equations are enforced:
\begin{subequations}
    \begin{align}
        \mathbf{G}(t_j) = \mathbf{P}^-(t_j) {\mathbf{H}^-(t_j)}^\top (\mathbf{H}^-(t_j) \mathbf{P}^-(t_j) {\mathbf{H}^-(t_j)}^\top \!\!+ \mathbf{R})^{-1}, \\
        \hat{\mathbf{x}}^+(t_j) = \hat{\mathbf{x}}^-(t_j) + \mathbf{G}(t_j) (\bar{\mathbf{y}}(t_j) - \mathbf{h}(\hat{\mathbf{x}}^-(t_j),t_j)),
        \label{eq:statePosterior}\\
        \mathbf{P}^+\!(t_j)\! =\! \left(\mathbf{I}\!-\!\mathbf{G}(t_j)\mathbf{H}^-(t_j)\right)\mathbf{P}^-(t_j)\left(\mathbf{I}\!-\!\mathbf{G}(t_j) \mathbf{H}^-(t_j)\right)^\top\!\! \!+\! \mathbf{G}(t_j) \mathbf{R} \mathbf{G}(t_j)^\top,
        \label{eq:covPosterior}
    \end{align}
\end{subequations}
where: $\mathbf{G}(t_j)$ is the Kalman gain matrix at $t_j$; $\mathbf{H}$ is the Jacobian of $\mathbf{h}$. The matrix $\mathbf{H}^-(t_j)$ is evaluated for $\mathbf{x}^-(t_j)$ as in the following:
\begin{equation}
    \mathbf{H}^-(t_j) =\frac{\partial \mathbf{h}}{\partial \mathbf{x}^{\top}}\bigg|_{\hat{\mathbf{x}}^-(t_j)}.
    \label{eq:obsEqsJac}
\end{equation}
When the KF is used to predict and correct the estimate of a set of variables describing the system dynamics, it is said that the procedure performs \textit{state estimation}. The procedure can be extended to perform \textit{joint estimation} as we will discuss in Sec. \ref{sec:jointEstimation}. As noted in the Introduction, the dimensionality of the observations ${\mathbf{y}}\in\mathbb{R}^{o}$ does not necessarily coincide with the dimensionality of the system state $\mathbf{x}\in\mathbb{R}^n$, thus potentially leading to $o \neq n $.

\subsection{State estimation of autonomous dynamical system}

In this subsection we present the main tools we combine with EKF to devise our proposed EKF-SINDy strategy, namely, a SINDy model and time-delay embedding. We thus explain how this combination can be exploited to perform the joint estimation of non-autonomous dynamical systems.

\subsubsection{Construction of the SINDy model}

The mapping $\mathbf{f}$ is only sometimes known, but, even if this is the case, the computation of $\mathbf{F}$ can be rather involved, while its numerical calculation can be computationally expensive. These reasons suggest to use SINDy to replace $\mathbf{f}$. Indeed, by expressing $\mathbf{f}$ as the linear combination of a set of predetermined functions collected in a library, as SINDy does, its (partial) derivatives can be simply obtained by combining the (partial) derivatives of the library functions. Analogously, even though $\mathbf{h}$ does not describe a dynamical system, it could be also identified using SINDy, thus greatly simplifying the computation of $\mathbf{H}$. However, in general the function $\mathbf{h}$ is directly set by means of more straightforward techniques. In the first case study proposed in this work, system state and observations represent the same dynamic quantities, and they do not require a distinct identification for $\mathbf{h}$. In the second case study, $\mathbf{h}$ has been determined from the time-delay embedding required by SINDy to (fully) learn $\mathbf{f}$.

For the sake of generality, we now consider using SINDy to model both $\mathbf{f}$ and $\mathbf{h}$. The training of SINDy is performed in an offline phase, prior to the online assimilation of incoming data performed through the EKF prediction-correction scheme. A detailed procedure is reported for $\mathbf{f}$, but the same strategy can also be applied to $\mathbf{h}$. The discussion is related to autonomous systems, i.e., systems without external forcing; in the next section, this assumption will be relaxed to allow for non-autonomous systems.

To apply the SINDy technique, it is first necessary to collect snapshots of the state vector to define the following matrix:
\begin{equation}
     \mathbf{X} = \begin{bmatrix}        
    \mathbf{x}^{\top}(t_1) \\
    \mathbf{x}^{\top}(t_2) \\
    \vdots \\
    \mathbf{x}^{\top}(t_T)
    \end{bmatrix} = 
    \begin{bmatrix}
        x_1(t_1) & x_2(t_1) & \cdots & x_n(t_1) \\
        x_1(t_2) & x_2(t_2) & \cdots & x_n(t_2) \\
        \vdots & \vdots & \ddots & \vdots \\
        x_1(t_T) & x_2(t_T) & \cdots & x_n(t_T)
    \end{bmatrix}.
    \label{eq:snapshotMatrix}
\end{equation}
where $x_i(t_j)$ is the i-th entry of $\mathbf{x}$ at the $j$-th time instant.

Similarly, a matrix $\dot{\mathbf{X}}$ is constructed by collecting the time derivatives $\dot{\mathbf{x}}$. 
If $\mathbf{f}$ is available, Eq. \eqref{eq:dynamicSystem} can be used to generate $\mathbf{X}$ and $\dot{\mathbf{X}}$. Otherwise, $\mathbf{f}$ can be identified from experimental data exploiting noise tolerant versions of SINDy, see, e.g., \cite{art:Messenger21,art:Fasel22,art:Hirsh22,art:Gao23, art:conti2024veni}.

A library $\boldsymbol{\Theta}_f(\mathbf{x})=\left[\theta_1(\mathbf{x}), \ldots, \theta_p(\mathbf{x})\right]\in\mathbb{R}^p$ of $p$ candidate functions to describe the dynamics
of the data is selected, and the matrix $\boldsymbol{\Theta}_f(\mathbf{X})\in\mathbb{R}^{T\times p}$ is constructed by applying $\boldsymbol{\Theta}_f$ to the rows of $\mathbf{X}$, e.g., as in the following:
\begin{equation}
\boldsymbol{\Theta}_f(\mathbf{X})=\begin{bmatrix}
        | & | & | & | & & | & | & \\
        1 & \mathbf{X} & \mathbf{X}^{P_2} & \mathbf{X}^{P_3} & \cdots & \text{sin}(\mathbf{X}) & \text{cos}(\mathbf{X}) & \cdots \\
        | & | & | & | & & | & | &
    \end{bmatrix}.
\label{eq:functionLabels}
\end{equation}

It is worth noting that any type of function can be utilised to form the function library. For instance, in Eq. \eqref{eq:functionLabels} we have employed constant, polynomial and trigonometric functions. Including polynomial functions is effective in identifying dominant dynamic behaviour, owing to the Taylor expansion of the function governing the system dynamics \cite{art:Champion19}. The way in which the quadratic nonlinearities $\mathbf{X}^{P_2}$ are expressed is now explicitly detailed:
\begin{equation}    \mathbf{X}^{P_2}=\begin{bmatrix}
        x^2_1(t_1) & x_1(t_1) x_2(t_1) & \cdots &  x^2_2(t_1) & \cdots &  x^2_n(t_1) \\
        x^2_1(t_2) &  x_1(t_2) x_2(t_2) & \cdots & x^2_2(t_2) & \cdots & x^2_n(t_2) \\
        \vdots & \vdots & \ddots & \vdots & \ddots & \vdots \\
        x^2_1(t_T) &  x_1(t_T) x_2(t_T) & \cdots & x^2_2(t_T) & \cdots & x^2_n(t_T)
    \end{bmatrix},
    \label{eq:quadraticNonLin}
\end{equation}
The definition of higher order nonlinearities can be done similarly.

The function $\mathbf{f}$ can be then expressed as the linear combination of these candidate functions. The weighting coefficients of the combination are stored in the matrix $\boldsymbol{\Xi}_f=\left[ \boldsymbol{\xi}^f_1, \boldsymbol{\xi}^f_2,\ldots, \boldsymbol{\xi}^f_n\right]$ with $\boldsymbol{\xi}^f_i\in \mathbb{R}^{p}$ for $i=1,\ldots,n$. In matrix form, such combination can be rewritten as:
\begin{equation}
    \dot{\mathbf{X}} = \boldsymbol{\Theta}_f(\mathbf{X}) \boldsymbol{\Xi}_f
    \label{eq:sindy},
\end{equation}
where the equality sign holds as $\boldsymbol{\Theta}_f$ admits an arbitrary large number of terms, thus allowing -- in principle -- to perfectly reconstruct the system dynamics. As only few terms of the function library are expected to be able to describe the dynamics of the system of interest, it is assumed that $\mathbf{f}$ admits a sparse representation in $\boldsymbol{\Theta}_f$. Consequently, a sparsity promoting regulariser $\mathcal{L}(\boldsymbol{\xi}'_i)$ is added in the least square regression, used to determine the weighting coefficients, according to:
\begin{equation}
\boldsymbol{\xi}^f_i=\underset{\boldsymbol{\xi}'_i}{\arg\min}\Vert \dot{\mathbf{X}}_i - \boldsymbol{\Theta}_f(\mathbf{X})\boldsymbol{\xi}'_i \Vert_2 + \mathcal{L}(\boldsymbol{\xi}'_i). 
\label{eq:weightSINDy}
\end{equation}

For instance, problem \eqref{eq:weightSINDy} could be solved via LASSO regression \cite{hastie2015statistical}, by considering $\mathcal{L}$ to be a sparsity-promoting $\ell_1-$norm, \i.e. $\mathcal{L}(\boldsymbol{\xi}'_i) = \delta_l \Vert \boldsymbol{\xi}'_i \Vert_1$, where $\delta_l$ is the parameter setting the strength of the LASSO regularisation. Alternatively, the Sequential Thresholded Least SQuares (STLSQ) algorithm \cite{art:Brunton16} could be used. It consists in iteratively solving \eqref{eq:weightSINDy} with ridge regression by using $\ell_2-$regularization, \i.e. $\Vert \boldsymbol{\xi}'_i \Vert_2$, regulated by a parameter $\delta_r$, and setting to zero the entries of $\boldsymbol{\xi}'_i$ with magnitude below a threshold $L$ at each iteration. 
As previously mentioned, SINDy can be used to identify both $\mathbf{f}$ and $\mathbf{h}$, even though $\mathbf{h}$ does not describe a system dynamics. To cope with that, a matrix $\mathbf{Y}$ collecting snapshots of $\mathbf{y}$ is constructed in addition to $\mathbf{X}$. Finally, Eq. \eqref{eq:sindy} is rewritten as:
\begin{equation}
    \mathbf{Y} = \boldsymbol{\Theta}_h(\mathbf{X}) \boldsymbol{\Xi}_h
    \label{eq:sindy_h},
\end{equation}
where the employed function library $\boldsymbol{\Theta}_h$ can be in general different from $\boldsymbol{\Theta}_f$. The weighting coefficient of the combination $\boldsymbol{\xi}^h_i$, stored in the matrix $\boldsymbol{\Xi}_h=\left[ \boldsymbol{\xi}^h_1, \boldsymbol{\xi}^h_2,\ldots, \boldsymbol{\xi}^h_o\right]$, are determined through a least square regression similar to what done in Eq. \eqref{eq:weightSINDy}, provided the substitution of $\dot{\mathbf{X}}_i$ by $\mathbf{Y}_i$, of $\boldsymbol{\Theta}_f$ by $\boldsymbol{\Theta}_h$.

Modelling $\mathbf{f}$ and $\mathbf{h}$ with SINDy enables to easily compute the Jacobian matrices $\mathbf{F}$ and $\mathbf{H}$, as it implies to calculate the derivatives of the functions of the library $\boldsymbol{\Theta}_f$ and $\boldsymbol{\Theta}_h$. Specifically, by combining Eqs. \eqref{eq:Jacobian_def} and \eqref{eq:obsEqsJac} with Eqs. \eqref{eq:sindy} and  \eqref{eq:sindy_h}, $\mathbf{F}$ and $\mathbf{H}$ assume the following forms:
\begin{subequations}
    \begin{align}
    \mathbf{F}^-(\tau) = \frac{\partial}{\partial \mathbf{x}^{\top}} \left[ \boldsymbol{\Xi}_f^{\top} \boldsymbol{\Theta}^{\top}_f(\mathbf{x})\right]_{\hat{\mathbf{x}}^-(\tau)}  = \boldsymbol{\Xi}^{\top}_f \frac{\partial \boldsymbol{\Theta}^{\top}_f(\mathbf{x})}{\partial \mathbf{x}^{\top}} \bigg|_{\hat{\mathbf{x}}^-(\tau)} ,
    \label{eq:evJacobianSINDy}\\
    \mathbf{H}^-(t_j) =\frac{\partial}{\partial \mathbf{x}^{\top}}\left[ \boldsymbol{\Xi}^{\top}_h \boldsymbol{\Theta}^{\top}_h(\mathbf{x}) \right]_{\hat{\mathbf{x}}^-(t_j)}= \boldsymbol{\Xi}^{\top}_h \frac{\partial \boldsymbol{\Theta}^{\top}_h(\mathbf{x})}{\partial \mathbf{x}^{\top}} \bigg|_{\hat{\mathbf{x}}^-(t_j)},
    \label{eq:obsJacobianSINDy}
\end{align}
\end{subequations}
marking how much the use of SINDy can simplify the application of the hybrid version of the EKF.

\subsubsection{Time-delay embedding}
\label{sec: timedelay}
\begin{figure}[t!]
    \centering
    \includegraphics[width=150mm]{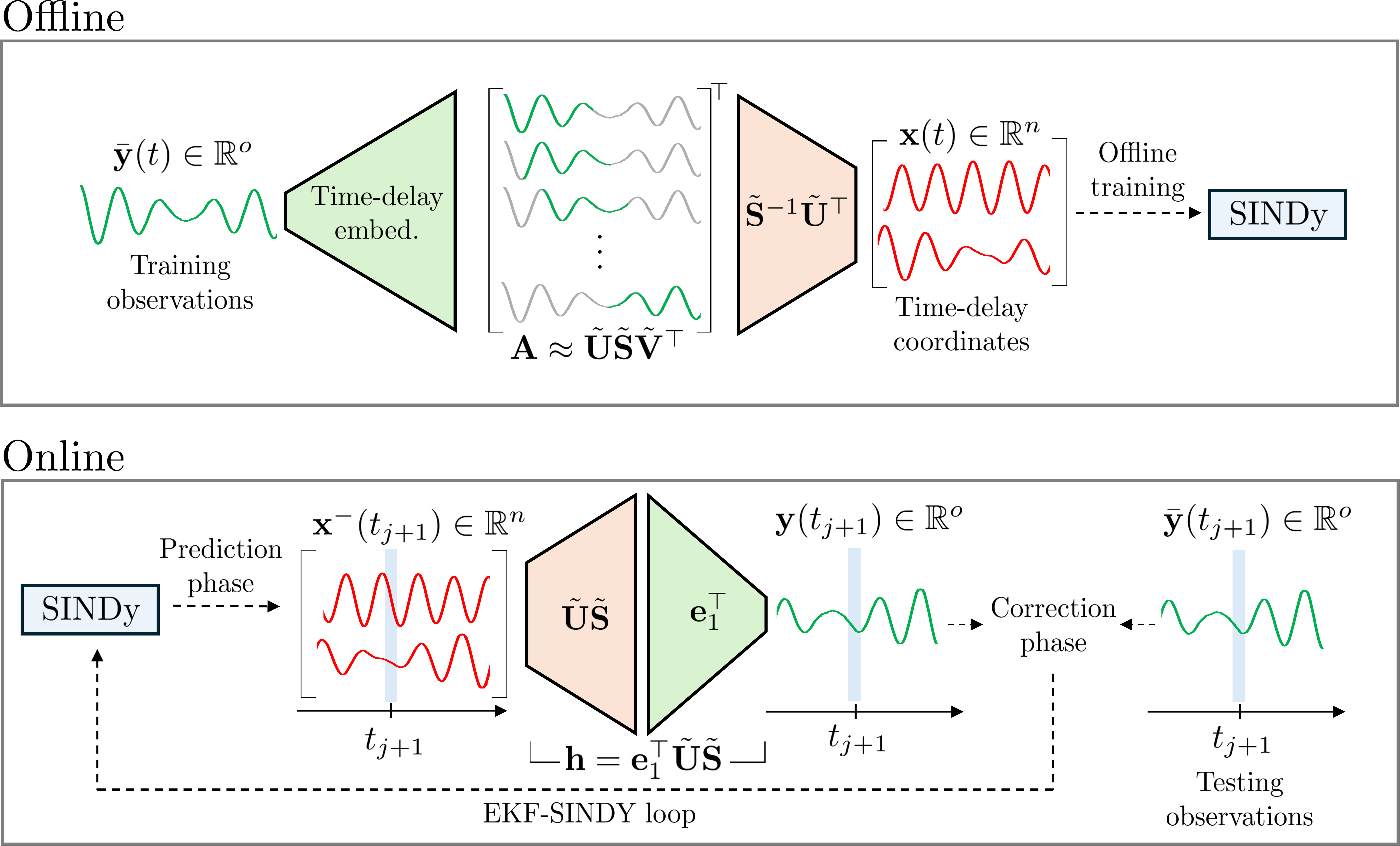}
    \caption{{Time-delay emebdding procedure. Over the offline phase, time-series observations of the system are time-delay embedded to form the Hankel matrix $\mathbf{A}$. The time-delay coordinates are then obtained by projecting the delayed signals onto the dominant $\eta$ left-singular values of the SVD decomposition (\i.e., columns of $\tilde{\mathbf{U}}$) and rescaled by dividing them by the corresponding singular values, stored in $\tilde{\mathbf{S}}$. Finally, SINDy is trained offline on these time-delay variables, approximating the full state space variables of the system under consideration (we stressed this aspect by indicating $\eta=n$, even though in general $\eta$ could be different from $n$). Once trained offline, SINDy could be employed in the online prediction phase to forecast the evolution of the time-delay variables, from which the observations can be reconstructed by back-projection (\i.e. by multiplication for $\tilde{\mathbf{U}}\tilde{\mathbf{S}}$). 
    Then the reconstructed observations can be simply extracted by looking at the first column, which contains the unshifted reconstructed signal, and compared to the actual testing observations (assimilated online) to perform the final correction phase. It follows that state-to-observation map, $\mathbf{h}$, is straightforwardly obtained as $\mathbf{h}=\mathbf{e}_1 ^\top\tilde{\mathbf{U}}\tilde{\mathbf{S}}$, where $\mathbf{e}_1$ is the basis vector extracting the first column of $\tilde{\mathbf{U}}\tilde{\mathbf{S}}$.}}
    \label{fig: timedelay}
\end{figure}

A possible drawback of SINDy is that its training requires knowledge of the full underlying state space variables. In many real-world applications, states could not be directly observed and some governing variables may be completely unobserved, resulting in observations ${\mathbf{y}}\in\mathbb{R}^{o}$ which do not match the system states ${\mathbf{x}}\in\mathbb{R}^{n}$. In the partial measurement scenario, with $o<n$, time-delay embedding techniques \cite{art:Brunton17} can be adopted to recover hidden state components that might be fundamental to describe the dynamics \cite{bakarji2023discovering}. Hereon we provide a description of this strategy -- see Fig.~\ref{fig: timedelay} for an overall sketch.
 
First, time-delay embedding method consists in lifting the low-dimensional observations $y(t)$ (considering $o=1$ for simplicity) into high-dimensional space by stacking time-shifted copies of the measurements $y(t)$ by means of vectors $\mathbf{a}(t;w,\zeta)=[y(t), y(t+\zeta),y(t+2\zeta),\ldots, y(t+(w-1)\zeta)]^\top\in\mathbb{R}^w$, where $w$ and $\zeta$ indicate the number of delay embeddings and the lag time, respectively. The resulting matrix $\mathbf{A}$ is termed as Hankel matrix and has the following form (taking $\zeta = 1$):
\begin{equation}
     \mathbf{A}=
     \begin{bmatrix}
        y(t_1) & y(t_2) & \cdots & y(t_{T-w}) \\
        y(t_2) & y(t_3) & \cdots & y(t_{T-w+1}) \\
        \vdots & \vdots & \ddots & \vdots \\
        y(t_w) & y(t_{w+1}) & \cdots & y(t_{T})
    \end{bmatrix} = 
         \begin{bmatrix}
        \mathbf{a}(t_1), \mathbf{a}(t_2), \ldots, \mathbf{a}(T-w)
    \end{bmatrix} \in \mathbb{R}^{w\times(T-w)}.
\label{eq:hankel}%
\end{equation}
Second, a Singular Value Decomposition (SVD) $\mathbf{A}=\mathbf{U}\mathbf{S}\mathbf{V}^\top$ is performed to extract time-delay coordinates approximating the Koopman operator in finite dimension \cite{art:Dylewsky22}, thereby providing a linearised description of the system dominant dynamics \cite{art:Champion19}. 
The matrices $\mathbf{U}\in\mathbb{R}^{w\times w}$ and $\mathbf{V} \in\mathbb{R}^{(T-w)\times (T-w)}$ are two orthogonal matrices whose columns are termed left and right singular vectors, respectively; ${\mathbf{S}}\in\mathbb{R}^{w\times (T-w)}$ is a pseudo-diagonal matrix collecting the associated singular values.
The reference system formed by 
the columns of $\mathbf{U}$ is diffeomorphic to the original attractor of the system dynamics, under the conditions given by the Takens'~embedding theorem \cite{proc:Takens81}. The number of time-delay coordinates often can be significantly reduced to $\eta\ll w$ and can be chosen according to the decay of the corresponding singular values, thus resulting in a truncated SVD: $\mathbf{A}\approx\tilde{\mathbf{U}}\tilde{\mathbf{S}}\tilde{\mathbf{V}}^\top$, with $\tilde{\mathbf{U}}\in\mathbb{R}^{w\times \eta},\tilde{\mathbf{S}}\in\mathbb{R}^{\eta\times \eta}, \tilde{\mathbf{V}}\in\mathbb{R}^{(T-w)\times \eta}$. 
The selected $\eta$ coordinates are expected to retrieve the unobserved, hidden, state variables, and therefore can be leveraged by SINDy in the offline training phase, to better capture the dynamics and to improve prediction performance \cite{bakarji2023discovering}. Ideally $\eta=n$, being $n$ the state dimensionality of the partially-observed system under consideration. 

An additional advantage of time-delay embedding is that it allows the reconstruction of the original variables through back-projection. This provides an explicit expression for the function $\mathbf{h}$ that maps the identified latent states to the corresponding observation vector and its Jacobian $\mathbf{H}$. This enables the online prediction phase to be performed in the time-delay embedding coordinates and to exploit low dimensional observations in the correction phase, thus maintaining the EKF-SINDy online procedure unchanged and avoiding the computational cost of computing the Jacobian \eqref{eq:obsEqsJac}.

The performance of the time-delay embedding could depend significantly on the parameters choices, such as the number of embeddings $w$ and the time lag $\zeta$. As general and qualitative guidance, time delay parameters should be chosen to allow for the attractor to properly unfold, thus occupying as much of the embedded phase space as possible \cite{bakarji2023discovering}. More quantitative techniques, based on relation or mutual information \cite{fraser1989information}, could be employed as in \cite{ma2006selection, kennel1992method}.

\subsection{Joint estimation of non-autonomous dynamical systems}
\label{sec:jointEstimation}

The functions $\mathbf{f}$ and $\mathbf{h}$ do not necessarily depend only on the variables $\mathbf{x}$ describing the system dynamics, but also on a set of parameters $\boldsymbol{\phi}\in\mathbb{R}^l$ having a precise physical meaning. For example, the response of a building to an earthquake excitation depends on the stiffness of its structural members whose value may change in time due to long-term degradation of concrete and/or extreme events.

The goal of joint estimation is to simultaneously estimate $\mathbf{x}$ and the parameters $\boldsymbol{\phi}$ \cite{art:Mariani05}. To reach the scope, we increase the number of state variables by augmenting the state vector as $\boldsymbol{\varkappa}=[\mathbf{x}^{\top},\boldsymbol{\phi}^{\top}]^{\top}$. We consider the case in which it is impossible to associate a dynamic description to the evolution of $\boldsymbol{\phi}$. In this setting, the prediction stage of $\boldsymbol{\phi}$ is performed by a random walk driven by a process noise $\mathbf{q}_{\phi}$ as in the following:
\begin{equation}
\dot{\boldsymbol{\phi}}=\mathbf{q}_{\phi},
 \label{eq:phiDyn}
\end{equation}
where $\mathbf{q}_{\phi}$ is assumed to be sampled from a white, zero mean stochastic process $\mathbf{q}_{\phi}\sim \mathcal{N}(\mathbf{0},\mathbf{Q}_{\phi})$ with a diagonal covariance matrix $\mathbf{Q}_{\phi}\in\mathbb{R}^{l\times l}$. Keeping the definition of $\mathbf{q}_x$ and $\mathbf{Q}_x$, we introduce a unique process noise vector $\mathbf{q}=[\mathbf{q}_x,\mathbf{q}_{\phi}]^{\top}$ and a unique diagonal covariance matrix $\mathbf{Q}$ collecting the entries of $\mathbf{Q}_x$ and $\mathbf{Q}_{\phi}$. Finally, we update the definition of $\mathbf{f}:\mathbb{R}^{n+l}\rightarrow\mathbb{R}^{n+l}$ to include both the function modelling the system dynamics and the parameter evolution ruled by Eq. \eqref{eq:phiDyn}.

In the correction stage, the mismatch (termed innovation) of the predictions with the incoming data is used to update both $\mathbf{x}$ and $\boldsymbol{\phi}$. The filter equations can be obtained by substituting $\mathbf{x}$ with $\boldsymbol{\varkappa}$. 

Also in case of joint estimation, it is possible -- and convenient -- to rely on SINDy. Specifically, the dependence of system dynamics on $\boldsymbol{\phi}$ can be accounted for by stacking snapshots $\boldsymbol{\varkappa}$ accounting for $I$ different realisations of $\boldsymbol{\phi}$, thereby obtaining a modified version $\mathbf{X}_{\varkappa}$ of the snapshot matrix. Thus, SINDy is employed to reconstruct the functional relation between $\mathbf{X}_{\varkappa}$ and $\dot{\mathbf{X}}$. It is unnecessary to modify the definition of $\dot{\mathbf{X}}$, as we have assumed to drive the dynamics of $\boldsymbol{\phi}$ with a random walk. 
It is worth noticing that also the definition function library $\boldsymbol{\Theta}_f(\boldsymbol{\varkappa})=[\theta_1(\boldsymbol{\varkappa}),\dots,\theta_p(\boldsymbol{\varkappa})]\in\mathbb{R}^p$ must be updated to account for the dependence on $\boldsymbol{\varkappa}$. Specifically, the matrix $\boldsymbol{\Theta}_f(\mathbf{X}_{\varkappa})$, defined in Eq. \eqref{eq:functionLabels}, is modified by applying $\boldsymbol{\Theta}$ to the rows of $\mathbf{X}_{\varkappa}$ obtaining:
\begin{equation}
\boldsymbol{\Theta}_f(\mathbf{X})=\begin{bmatrix}
        | & | & | & | & & | & | & \\
        1 & \mathbf{X}_{\varkappa} & \mathbf{X}_{\varkappa}^{P_2} & \mathbf{X}_{\varkappa}^{P_3} & \cdots & \text{sin}(\mathbf{X}_{\varkappa}) & \text{cos}(\mathbf{X}_{\varkappa}) & \cdots \\
        | & | & | & | & & | & | &
    \end{bmatrix}.
\label{eq:functionLabels_phi}
\end{equation}
Looking at the quadratic nonlinearities, Eq. \eqref{eq:quadraticNonLin} is updated to include the dependence on $\boldsymbol{\phi}$ as follows:
\begin{equation}    \mathbf{X}^{P_2}_{\varkappa}=\begin{bmatrix}
        x^2_1(t_1) & x_1(t_1) x_2(t_1) & \cdots & x_1(t_1)\phi_1(t_1) & \cdots &  x^2_n(t_1) & \cdots & {\phi}^2_l(t_1) \\
        x^2_1(t_2) & x_1(t_2) x_2(t_2) & \cdots & x_1(t_2)\phi_1(t_2) & \cdots &  x^2_n(t_2) & \cdots & {\phi}^2_l(t_2) \\
        \vdots & \vdots & \ddots & \vdots & \ddots & \vdots & \ddots & \vdots \\
         x^2_1(t_T) & x_1(t_T) x_2(t_T) & \cdots & x_1(t_T)\phi_1(t_T) & \cdots &  x^2_n(t_T) & \cdots & {\phi}^2_l(t_T)
    \end{bmatrix}.
    \label{eq:quadraticNonLin_phi}
\end{equation}
Higher order linearities and trigonometric functions can be similarly treated.

Notably, setting $\dot{\boldsymbol{\phi}}=0$ implies that the parameter values can be updated only in the correction phase of the EKF. This is allowed by updating Eq. \eqref{eq:evJacobianSINDy} to define the Jacobian $\mathbf{F}\in\mathbb{R}^{(n+l)\times(n+l)}$ including the dependence of system dynamics on $\boldsymbol{\phi}$ as in the following:
\begin{equation}
    \mathbf{F}^-(\tau) =
\begin{pNiceMatrix}[]
& \boldsymbol{\Xi}_f^{\top} \frac{\partial \boldsymbol{\Theta}^{\top}_f(\boldsymbol{\varkappa})}{\partial \mathbf{x}^{\top}}\bigg|_{\hat{\boldsymbol{\varkappa}}^-(\tau)} & & & \boldsymbol{\Xi}_f^{\top} \frac{\partial 
\boldsymbol{\Theta}^{\top}_f(\boldsymbol{\varkappa})}{\partial \boldsymbol{\phi}^{\top}}\bigg|_{\hat{\boldsymbol{\varkappa}}^-(\tau)} &
\\
& \mathbf{0} & &  & \mathbf{0} &\\
\CodeAfter
\UnderBrace[shorten,yshift=1.5mm]{last-1}{last-3}{\partial / \partial \mathbf{x}}
\UnderBrace[shorten,yshift=1.5mm]{last-4}{last-last}{\partial / \partial \boldsymbol{\phi}}
\end{pNiceMatrix},
\vspace{6mm}
\end{equation}
\noindent where the rows of $\mathbf{F}$ containing $\mathbf{0}$ entries, with $\mathbf{0}$ representing null matrices of suitable dimensions, are associated with the random walk equations. 
When used to model $\mathbf{h}$ for joint estimation, SINDy instead reconstructs the mapping between $\mathbf{X}_{\varkappa}$ and $\mathbf{Y}$. Also the definition of $\mathbf{Y}$ must not be modified as system parameters can not, in general, be observed. We highlight that the procedure can not explicitly handle the dependence on other parameters not included in $\boldsymbol{\phi}$, for which it is impossible to provide a separate estimate. This issue will be addressed in a future extension of the current work.

Up to now, we have considered autonomous system. However, the procedure can be easily extended to deal with non-autonomous systems by allowing for the presence of an external forcing term $\mathbf{b}\in\mathbb{R}^n$. Eqs. \eqref{eq:dynamicSystem} and \eqref{eq:obsEqDiscrete} modify in the following way:
\begin{subequations}
    \begin{align}
    \dot{\mathbf{x}}(t) &= \mathbf{f}(\mathbf{x}(t),\mathbf{b}(t)),
    \label{eq:dynamicSystem_b}\\
    \mathbf{y}(t) &= \mathbf{h}(\mathbf{x}(t),\mathbf{b}(t)).
    \label{eq:obs_noise_b}
\end{align}
\end{subequations}

Similar modifications apply to Eqs. \eqref{eq:dynamicSystem_noise} and \eqref{eq:obs_noise}. The other equations of the filter remain unchanged if we assume to know $\mathbf{b}$ that is in the case of an exogenous input. This is true, for instance, if we consider a seismic excitation thanks to the presence of monitoring networks \cite{art:DAlessandro19}. If $\mathbf{b}$ was unknown, it would be possible to modify the formulation of the filter to perform input--state--parameter estimation \cite{art:EftekharAzam15,art:Dertimanis19,art:Castiglione20,art:EbrahimzadehHassanabadi23}. Also this further extension of the method will be considered in future.

The external forcing term can be handled via SINDy by constructing a snapshot matrix $\mathbf{X}_{\varkappa b}$ assembling different snapshots $[\boldsymbol{\varkappa}(t_j),\mathbf{b}(t_j)]^T$ determined for several realisations of $\boldsymbol{\phi}$. As a result, SINDy is employed to reconstruct the functional relation between $\mathbf{X}_{\varkappa b}$ and $\dot{\mathbf{X}}$ (and between $\mathbf{X}_{\varkappa b}$ and $\mathbf{Y}$ when used to model $\mathbf{h}$). Once again, the function library $\boldsymbol{\Theta}_f$ should be updated to account for the dependence on $\mathbf{b}$, similarly to what has been previously done in Eqs. \eqref{eq:functionLabels} and \eqref{eq:quadraticNonLin} to consider the dependence on $\boldsymbol{\phi}$.

The complete procedure for the joint estimation of non-autonomous dynamical system are reported in Algorithms \ref{al:offPhase} and \ref{al:onPhase}. In the algorithms, we consider the most general case in which both $\mathbf{f}$ and $\mathbf{h}$ are identified by using SINDy. The training of the SINDy model, and therefore the assembling of the snapshot matrices, must be done in an offline phase a priori with respect to the application of the procedure. Algorithm \ref{al:offPhase} is devoted to the illustration of this phase assuming to know $\mathbf{f}$ and $\mathbf{h}$. This case is of interest since it will be considered in the results section. As previously remarked, it is also possible to estimate $\mathbf{f}$ and $\mathbf{h}$ using noise tolerant versions of SINDy. During the online phase, the procedure is applied to update the estimates of the system by using incoming data. Algorithm \ref{al:onPhase} refers to this second part.

\begin{algorithm}[t!]
\hspace*{\algorithmicindent} \textbf{Input}: libraries of candidate functions $\boldsymbol{\Theta}_f([\boldsymbol{\varkappa},\mathbf{b}]^{\top})$ and $\boldsymbol{\Theta}_h([\boldsymbol{\varkappa},\mathbf{b}]^{\top})$ \\
\hspace*{\algorithmicindent}
\textbf{Output}: SINDy models $\mathbf{f}\approx\boldsymbol{\Xi}^{\top}_f\boldsymbol{\Theta}^{\top}_f([\boldsymbol{\varkappa},\mathbf{b}]^{\top})$ and $\mathbf{h}\approx \boldsymbol{\Xi}^{\top}_h\boldsymbol{\Theta}^{\top}_h([\boldsymbol{\varkappa},\mathbf{b}]^{\top})$
\begin{algorithmic}[1]
\For{$\iota=1,\ldots,I$}
\State Sample $\boldsymbol{\phi}_{\iota}$
\For{$t_j=t_1,\ldots,t_T$}
\State Evolve $ \dot{\mathbf{x}}\!\! =\!\! \mathbf{f}(\boldsymbol{\varkappa},\mathbf{b})$  from $t_{j\!-\!1}$ to $t_j$, with $\boldsymbol{\varkappa}(t_j)\!\!=\!\![\mathbf{x}(t_j),\boldsymbol{\phi}_{\iota}]^{\top}\!\!$ and $\dot{\boldsymbol{\phi}}_{\iota}\!\!=\!0\!$
\State Collect $[\boldsymbol{\varkappa}(t_j),\mathbf{b}(t_j)]^{\top}$ in $\mathbf{X}_{\varkappa b}$ 
\State Collect $\dot{\mathbf{x}}(t_j)$ in $\dot{\mathbf{X}}$
\State Compute $\mathbf{y}(t_j) = \mathbf{h}(\mathbf{x}(t_j),\mathbf{b}(t_j))$
\State Collect $\mathbf{y}$ in $\mathbf{Y}$
\EndFor
\EndFor
\State Set $\boldsymbol{\Theta}_f(\mathbf{X}_{\varkappa b})$ and $\boldsymbol{\Theta}_h(\mathbf{X}_{\varkappa b})$
\State Compute $\boldsymbol{\xi}^f_i=\underset{\boldsymbol{\xi}'_i}{\arg\min}\Vert \dot{\mathbf{X}}_i - \boldsymbol{\Theta}_f(\mathbf{X}_{\varkappa b})\boldsymbol{\xi}'_i \Vert_2 + \mathcal{L}(\boldsymbol{\xi}'_{i})$ for $i=1,\ldots,n$
\State Compute $\boldsymbol{\xi}^h_{i}=\underset{\boldsymbol{\xi}'_{i}}{\arg\min}\Vert {\mathbf{Y}}_{i} - \boldsymbol{\Theta}_{h}(\mathbf{X}_{\varkappa b})\boldsymbol{\xi}'_i \Vert_2 + \mathcal{L}(\boldsymbol{\xi}'_{i})$ for $i=1,\ldots,o$
\State Assemble $\boldsymbol{\Xi}_f=[\boldsymbol{\xi}^f_i]$ for $i=1,\ldots,n$
\State Assemble $\boldsymbol{\Xi}_h=[\boldsymbol{\xi}^h_i]$ for $i=1,\ldots,o$
\end{algorithmic}
\caption{EKF-SINDy joint estimation ($\boldsymbol{\varkappa}=[\mathbf{x},\boldsymbol{\phi}]^{\top}$) of non-autonomous dynamical systems: offline phase ($\mathbf{f}$ and $\mathbf{h}$ available)}
\label{al:offPhase}
\end{algorithm}

\begin{algorithm}[t!]
\hspace*{\algorithmicindent} \textbf{Input}: SINDy models $\mathbf{f}\approx\boldsymbol{\Xi}^{\top}_f\boldsymbol{\Theta}^{\top}_f([\boldsymbol{\varkappa},\mathbf{b}]^{\top})$ and $\mathbf{h}\approx\boldsymbol{\Xi}^{\top}_h \boldsymbol{\Theta}^{\top}_h([\boldsymbol{\varkappa},\mathbf{b}]^{\top})$; sequential measurements $\bar{\mathbf{y}}(t_{1}),\ldots,\bar{\mathbf{y}}(t_{T})$ \\
\hspace*{\algorithmicindent} \textbf{Output}: Estimates of $\hat{\boldsymbol{\varkappa}}^+(t_1),\ldots,\hat{\boldsymbol{\varkappa}}^+(t_T)$, and $\hat{\mathbf{P}}^+(t_1),\ldots,\hat{\mathbf{P}}^+(t_T)$ accounting for incoming data
\begin{algorithmic}[1]
\State Set $\mathbf{Q}$ and $\mathbf{R}$
\For{$t_j=t_1,\ldots,t_T$}
\[
\text{Predictor phase}
\]
\State $\mathbf{F}^+(t_{j\!-\!1}) =
\begin{bmatrix}
\boldsymbol{\Xi}^{\top}_f\frac{\partial \boldsymbol{\Theta}^{\top}_f([\boldsymbol{\varkappa},\mathbf{b}]^{\top})}{\partial \mathbf{x}^{\top}} \bigg|_{[\hat{\boldsymbol{\varkappa}}^+,\mathbf{b}](t_{j\!-\!1})} & \boldsymbol{\Xi}^{\top}_f\frac{\partial \boldsymbol{\Theta}^{\top}_f([\boldsymbol{\varkappa},\mathbf{b}]^{\top})}{\partial \boldsymbol{\phi}^{\top}} \bigg|_{[\hat{\boldsymbol{\varkappa}}^+,\mathbf{b}](t_{j\!-\!1})} \\
\mathbf{0} & \mathbf{0}
\end{bmatrix}$
\State $\hat{\boldsymbol{\varkappa}}^-(t_j) = \hat{\boldsymbol{\varkappa}}^+(t_{j\!-\!1})+\Delta t ~\boldsymbol{\Xi}^{\top}_f\boldsymbol{\Theta}^{\top}_f([\hat{\boldsymbol{\varkappa}}^+,\mathbf{b}]^{\top}(t_{j\!-\!1})) $
\State $\mathbf{P}^-\!(t_j)\! = \!\mathbf{P}^+\!(t_{j\!-\!1})+\Delta t\left(\mathbf{F}^+\!(t_{j\!-\!1})\mathbf{P}^+\!(t_{j\!-\!1})\! + \!\mathbf{P}^+\!(t_{j\!-\!1}) \mathbf{F}^{+^\top}\!\!(t_{j\!-\!1})\! +\! \mathbf{Q}\right)$
\[
\text{Corrector phase}
\]
\State $\mathbf{H}^-(t_j) =\begin{bmatrix}
\boldsymbol{\Xi}^{\top}_h\frac{\partial \boldsymbol{\Theta}^{\top}_h([\boldsymbol{\varkappa},\mathbf{b}]^{\top})}{\partial \mathbf{x}^{\top}} \bigg|_{[\hat{\boldsymbol{\varkappa}}^-,\mathbf{b}](t_j)} &  \boldsymbol{\Xi}^{\top}_h\frac{\partial \boldsymbol{\Theta}^{\top}_h([\boldsymbol{\varkappa},\mathbf{b}]^{\top})}{\partial \boldsymbol{\phi}^{\top}} \bigg|_{[\hat{\boldsymbol{\varkappa}}^-,\mathbf{b}](t_j)}
\end{bmatrix}$
\State $\mathbf{G}(t_j) = \mathbf{P}^-(t_j) {\mathbf{H}^-(t_j)}^\top (\mathbf{H}^-(t_j) \mathbf{P}^-(t_j) {\mathbf{H}^-(t_j)}^\top \!\!+ \mathbf{R})^{-1}$
\State $\hat{\boldsymbol{\varkappa}}^+(t_j) = \hat{\boldsymbol{\varkappa}}^-(t_j) + \mathbf{G}(t_j) (\bar{\mathbf{y}}(t_j) - \boldsymbol{\Xi}^{\top}_h \boldsymbol{\Theta}^{\top}_h([\hat{\boldsymbol{\varkappa}}^-,\mathbf{b}]^{\top}(t_j)))$
\State $\mathbf{P}^+(t_j)\! =\! \left(\mathbf{I}-\mathbf{G}(t_j)\mathbf{H}^-\!(t_j)\right)\mathbf{P}^-(t_j)\left(\mathbf{I}-\mathbf{G}(t_j) \mathbf{H}^-\!(t_j)\right)^\top \!+ \mathbf{G}(t_j) \mathbf{R} \mathbf{G}(t_j)^\top$
\EndFor
\end{algorithmic}
\caption{EKF-SINDy joint estimation ($\boldsymbol{\varkappa}=[\mathbf{x},\boldsymbol{\phi}]^{\top}$) of non-autonomous dynamical systems: online phase}
\label{al:onPhase}
\end{algorithm}

\section{Numerical Results}
\label{sec:results}

\subsection{Shear building under seismic excitations} \label{sec:shearBuildingResults}

\begin{figure}[b]
\centering
\includegraphics[width=140mm]{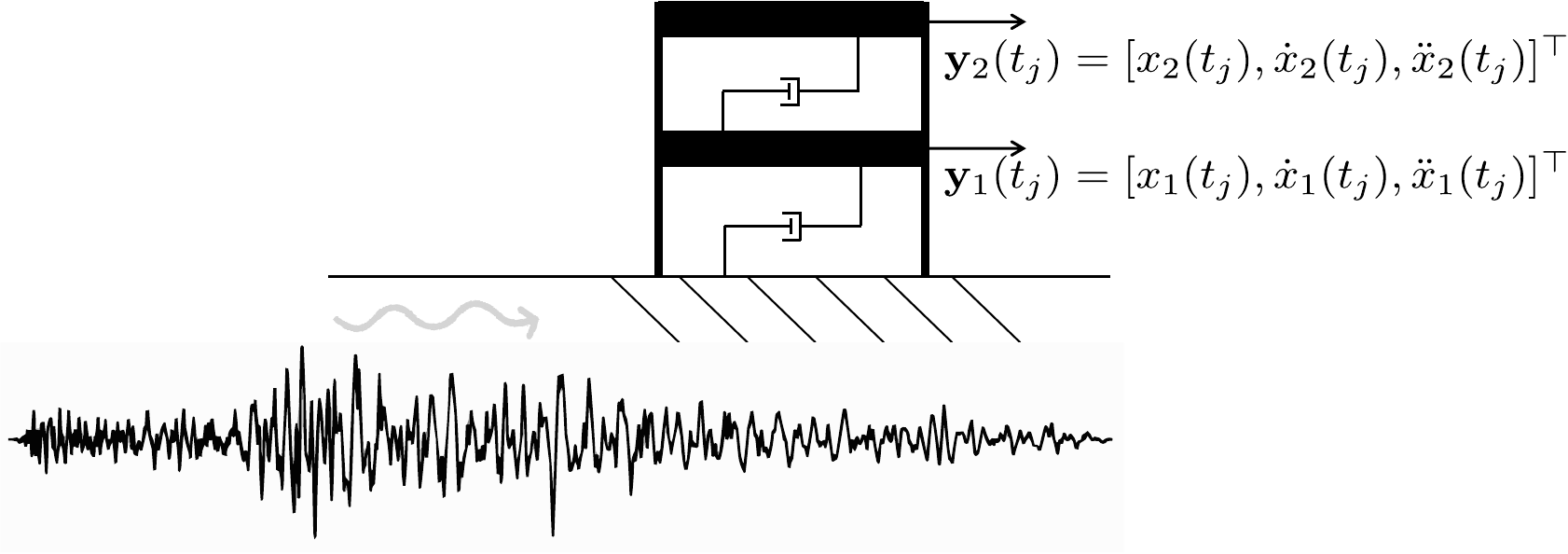}
\caption{\footnotesize Shear building model. Horizontal displacements, velocities, and accelerations at the two floors are recorded.}
\label{fig:shearModel}
\end{figure}

A first numerical case deals with a shear building under seismic excitations. Shear building models are effective in describing the response to lateral excitations, like seismic loads, whenever floors have a sufficiently high out-of-plane stiffness. If this requirement is met, e.g. by way of a minimum slab thickness, the adoption of shear building models in the design phase is allowed by standards such as Eurocode 8 \cite{code:EC8}. Here, we have examined the response to seismic excitation of a $2$ storey building assuming regularity in the floor mass and stiffness distribution, such that torsional effects can be neglected and a decoupling of the building response to lateral excitation along the two in-plan directions can be exploited. We thus have employed one degree-of-freedom (dof) per floor to model the response of the building along the horizontal direction. We have assumed the monitoring system to be deployed to record storey displacements $x_i(t)$, velocities $\dot{x}_i(t)$ and accelerations $\ddot{x}_i(t)$, with $i=1,2$, at discrete time steps, obtaining a $6$-dimensional observation vector $\mathbf{y}$. A schematic representation of the shear building model is reported in Fig. \ref{fig:shearModel}.

The dynamic response of the structure is described by the following equations:
\begin{subequations}
    \begin{align}
    m\ddot{x}_1(t) + c_1\dot{x}_1(t) + k (2x_1(t) - x_2(t)) = -m b(t) 
    \\
    m\ddot{x}_2(t) + c_2\dot{x}_2(t) + k (x_2(t) - x_1(t)) = -m b(t) 
\end{align}
\label{eq:shearBuilding}
\end{subequations}
where: $m=625$ ton is the mass of each floor; $k$ is the inter-storey stiffness; $b(t)$ is the ground acceleration signal; $c_1$ and $c_2$ are coefficients set to have $1\%$ damping over the two eigenmodes of the building model.

Seismic signals have been taken from the STEAD database \cite{art:STEAD}. Each seismogram $b(t)$ features the same time duration $60$~s and a sampling time step of $ 0.01$ s. Here, we have subsampled each signal getting a final sampling of $ 0.001$ s to avoid stability issue with the Euler forward integration. The dynamical system is not autonomous due to the presence of this external forcing term. We assume to know the input thanks to the seismic monitoring networks present in many countries, such as in Italy \cite{art:Pierleoni18,art:DAlessandro19}. On the contrary, we have made the assumption that the storey stiffness is unknown. For instance, when considering reinforced concrete structures, it is challenging to predict how cracking affects the moment of the inertia of column sections, or whether the modelling choices for boundaries have been adequate. Similarly, it is hard to guess how long term degradation effects may have altered the building properties over time. For this reason, we treat the value of the inter-storey stiffness as a stochastic variable featuring a uniform probability density function between $0.5\cdot 10^6$ kN$/$m and $2\cdot 10^6$ kN$/$m. Hence, the stiffness $k$ is the parameter to be here estimated together with the state vector $\mathbf{x}=[x_1, x_2, \dot{x}_1, \dot{x}_2]^{\top}$. The second order Eq. \eqref{eq:shearBuilding} can be rewritten at first order using these coordinates, and the dynamics of the system is accordingly described by $n=4$ variables corresponding to the lateral storey displacements and velocities.

To perform the joint estimation of this system, we have constructed a SINDy model for the function $\mathbf{f}$ describing at first order the system,  as detailed in Algorithm \ref{al:offPhase}, obtaining $\mathbf{f}\approx\boldsymbol{\Xi}^{\top}_f\boldsymbol{\Theta}^{\top}_f([\boldsymbol{\varkappa},b]^{\top})$, with $\boldsymbol{\varkappa}=[\mathbf{x},k]^{\top}$. Specifically, the snapshot matrix $\mathbf{X}_{\varkappa b}$ has been constructed by collecting in time the values of the augmented state vector $[\mathbf{x},k]^{\top}$ and for $I=20$ different values of $k$. Polynomial terms up to the second order have been employed in the definition of the function library $\boldsymbol{\Theta}_f$, as it can be demonstrated that they can describe the system dynamics, see, e.g. \cite{art:RUENG24}. Hence, we have obtained a model for the observation operator $\mathbf{h}$ as in the following
\begin{equation}
    \mathbf{y}=
    \begin{bmatrix}
         \mathbf{y}^x \\
         \mathbf{y}^{\dot{x}} \\
         \mathbf{y}^{\ddot{x}}
    \end{bmatrix} =
    \begin{bmatrix}
         \mathbf{B}_x\mathbf{x} \\
         \mathbf{B}_{\dot{x}}\dot{\mathbf{x}} \\
         \mathbf{B}_{\ddot{x}}(\boldsymbol{\Xi}^{\top}_f\boldsymbol{\Theta}^{\top}_f([\boldsymbol{\varkappa}^{\top},b]^\top))
    \end{bmatrix},
\label{eq:obsEqShearBuild}
\end{equation}
\noindent
where: $\mathbf{B}_x$ and $\mathbf{B}_{\dot{x}}$ are the Boolean matrices that relate displacements and velocities to the observations of the system; $\mathbf{B}_{\ddot{x}}$ is the Boolean matrix that relates the accelerations computed by SINDy to the observed accelerations. The computation of the Jacobian matrices $\mathbf{F}$ and $\mathbf{H}$ has been done as described in Eqs. \eqref{eq:evJacobianSINDy} and \eqref{eq:obsJacobianSINDy} taking advantage of the use of polynomial functions in the library $\boldsymbol{\Theta}_f$.

The least square regression used to determine the SINDy weighting coefficients, see Eq. \eqref{eq:weightSINDy}, is managed through the python open source package \texttt{PySINDy} \cite{art:Kaptanoglu22}. Sparsity is promoted using STLSQ with threshold $L=10^{-2}$ and regularisation strength $\delta_r=0.05$. The value for $\delta_r$ is the default setting in the PySINDy package. In the considered case studies, we have found that adjusting $L$ alone has been sufficient to control the sparsity of the system.
As expected, the adopted library of second-order polynomials is perfectly suitable to model the dynamics of the system, accounting also for the dependence on $k$.

Hence, we have taken advantage of the combination of SINDy with EKF to perform the joint estimation of shear building whose stiffness is unknown. Specifically, we have taken as initial guess $k=1.01\cdot 10^6$ kN$/$m, overestimating by $20\%$ the target value. As previously mentioned, displacements, velocities and accelerations of both floors have been recorded. Although an identifiability assessment has not been carried out, see e.g. \cite{art:Chatzis15}, the full observation of the system guarantees the possibility to estimate $k$. In future works, the intention is to explicitly perform an identifiability assessment, for example through the software DAISY \cite{art:Bellu07}.

To simulate the outcome of a real monitoring system, signals have been corrupted with white noise featuring a signal to noise ratio equal to $15$. Such level and type of noise is compatible with the use of, e.g., MEMS sensors \cite{proc:DAlessandro17}. Before running the online stage of Algorithm \ref{al:onPhase}, we have initialised the state covariance matrix $\mathbf{P}(0)$, the process noise covariance matrix $\mathbf{Q}$, and the observation noise covariance matrix $\mathbf{R}$ as diagonal matrices. The values of the diagonal entries of these matrices are specified in Appendix B. The tuning of these quantities has been carried out through a trial-and-error procedure (a practical guidance to expedite this procedure is included in Appendix B as well); however, automatic tuning procedures based on Bayesian Optimisation with Gaussian Processes \cite{art:Chen24}, on genetic algorithms \cite{art:Rapp03} or on swarm intelligence \cite{art:Laamari15}, can be utilised.

\begin{figure}[b!]
\centering
\includegraphics[width=145mm]{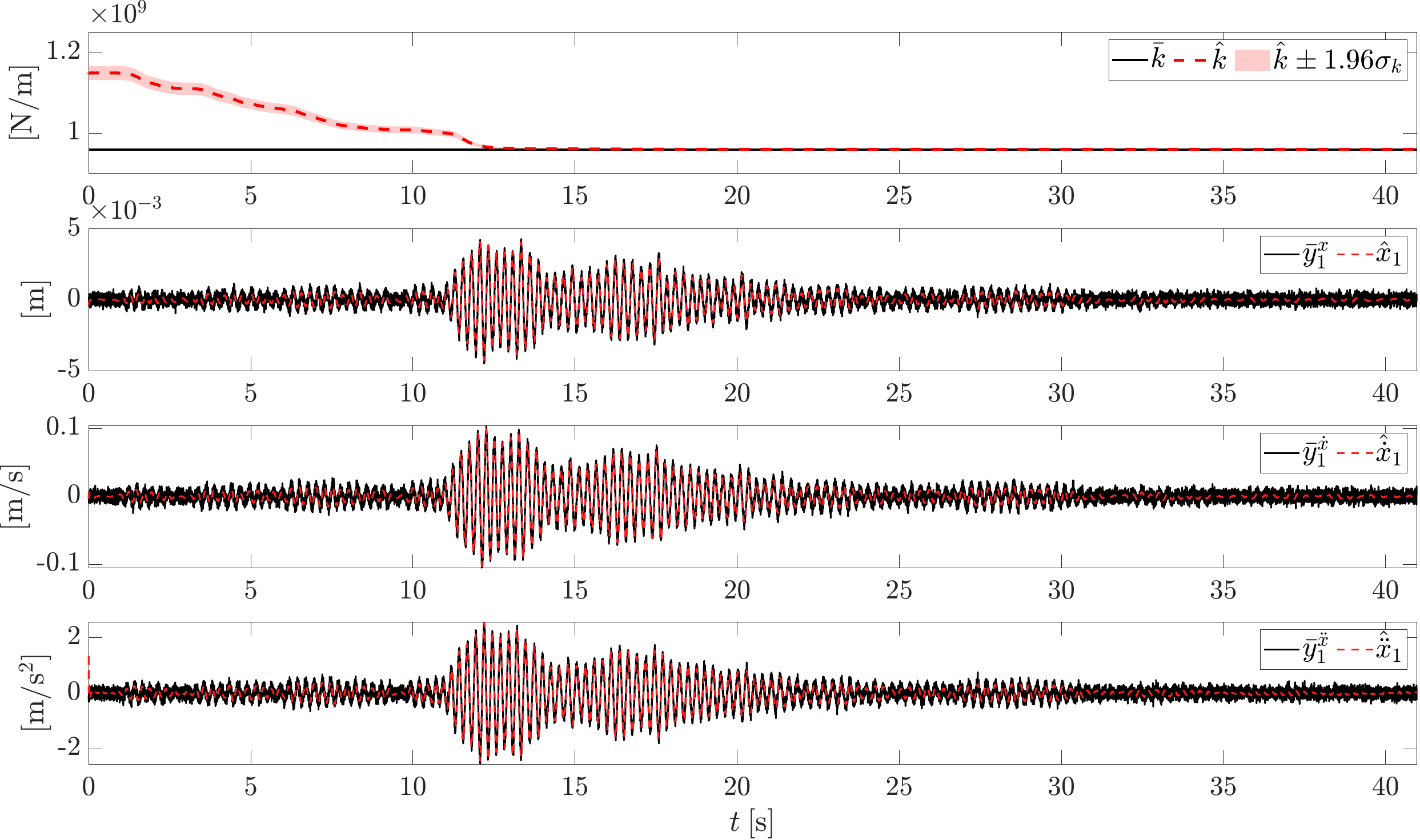}
{\caption{\footnotesize 
Shear building. In the top graph, the estimated parameter $\hat{k}$ is plotted in red against the target $\bar{k}$ in black. The red shaded area represents $95\%$ confidence interval of the estimates, determined using the posterior covariance. In the other graphs, the evolution of the system response tracked by the filter is plotted in red against the real system dynamics in black. For the sake of presentation, the outcomes of the $60$~s analyses are truncated after $40$ s.}\label{fig:2_dof_results}}
\end{figure}

\begin{figure}[b!]
\centering
\includegraphics[width=145mm]{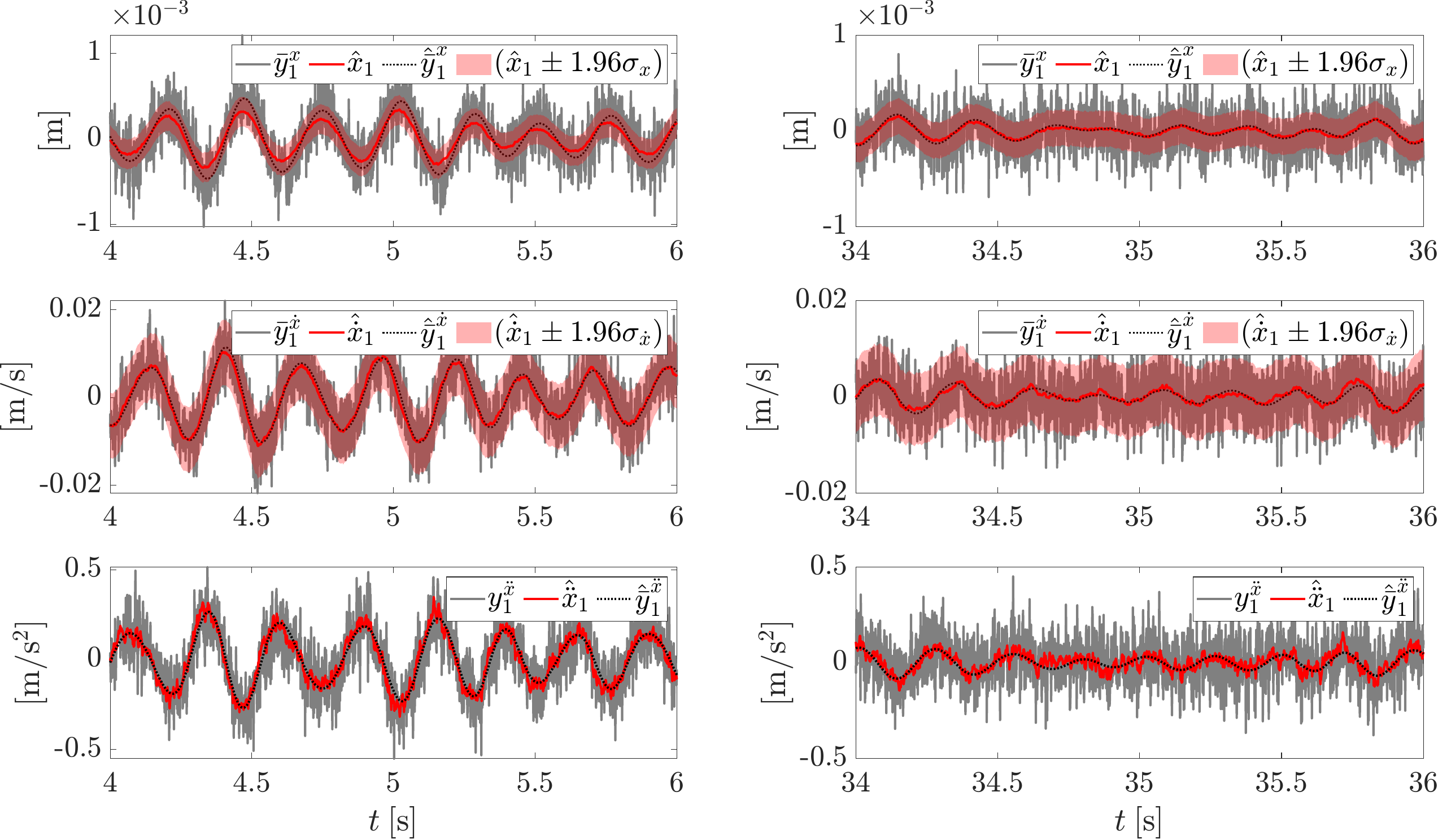}
{\caption{\footnotesize 
Shear building. Zoom in of Fig. \ref{fig:2_dof_results}: comparison between the evolution in time of the system response tracked by the filter (red dotted lines) and target one (black lines) for incorrect (on the left) and updated (on the right) value of $k$. Red shaded areas represent the $95\%$ confidence intervals determined using the posterior covariance. Black solid lines represent the noise corrupted signals, while black dotted lines their uncorrupted version.}\label{fig:2_dof_results_zoom_in}}
\end{figure}

The outcome of the estimation procedure is summarised in Fig. \ref{fig:2_dof_results} for a seismogram randomly picked from the STEAD database. The mean value of the estimated parameter $\hat{k}$ converges to the target $\bar{k}$ in roughly 20~s reducing the uncertainty of the estimation, measured by the standard deviation $\sigma_k$ determined using the posterior covariance computed by the EKF. The system dynamics is precisely tracked as well, as illustrated in Fig. \ref{fig:2_dof_results_zoom_in} in terms of a close-up of the kinematic quantities. A confidence interval for the acceleration estimates is not provided, because accelerations are not included in the state vector $\mathbf{x}$. The procedure is strongly tolerant to noise, being the filter estimates (for example $\hat{x}_1$) much closer to the noise-free version of the signals ($\hat{\bar{y}}^{x}_1$) than the noise corrupted versions ($\bar{y}^{x}_1$), despite that in the first part of the analysis $k$ is not correctly calibrated yet. Clearly, the matching improved even more after the first $20$ s.

Beyond the excellent estimation capabilities, the most intriguing aspects of the proposed procedure are its ease of use and the short execution time. Specifically, the Jacobian matrices have been straightforwardly determined (the reader may refer to \cite{art:RUENG24} to understand the level of effort required to compute the Jacobian matrices for a simple system like this). Regarding the computational aspects, the execution of the procedure is approximately $40$ times faster than the physical process (lasting $1.5$~s when run on a workstation equipped with an Intel (R) Core$^{\text{TM}}$, i7-2600 CPU @ 3.4 GHz with 16 GB RAM), moving towards real time monitoring in civil and mechanical engineering applications.

\subsection{Partially observed nonlinear dynamical system} \label{sec:resonatorResults}

A second test case deals with a partially observed nonlinear dynamical system. 
The system comprises two coupled oscillators characterised by linear stiffness coefficient $k_i$, and linear damping coefficient $c_i$, with $i=1,2$. These oscillators are linked via linear and quadratic coupling terms controlled by the coefficients $\alpha$ and $\beta$, respectively. The coefficient $\gamma$ dictates instead the contribution of a cubic nonlinearity to the motion of one oscillator. The equations governing the free vibrations of the system are:
\begin{subequations}
    \begin{align}
    m_1\ddot{z}_1(t) + c_1 \dot{z}_1(t) + k_1 z_1(t) + \alpha z_2(t) = 0,    \label{eq:coupledOscillator1}\\
    m_2\ddot{z}_2(t) + c_2 \dot{z}_2(t) + k_2 z_2(t) + \gamma z^3_2(t) + \alpha z_1(t) + \beta z^2_1(t)  = 0,\label{eq:coupledOscillator2}
\end{align}
\label{eq:coupledOscillator12}
\end{subequations}

\noindent where: $m_1=1$; $m_2=1$; $k_1=1.0$; 
$c_1=2\cdot 10^{-2}$; $c_2=1.95\cdot 10^{-2}$; $\alpha=-10^{-1}$; $\beta=2 \cdot10^{-3}$; $\gamma=10^{-3}$. The time axis, of length $200.0$, is discretised with a time step of $\Delta t = 10^{-2}$. Both the parameters and the time axis have been made dimensionless. Similarly to what considered in Sec. \ref{sec:shearBuildingResults}, the second order Eqs. \eqref{eq:coupledOscillator12} have been rewritten at first order with a simple change of coordinate, \i.e. $(\tilde{z}_1, \tilde{z}_2, \tilde{z}_3, \tilde{z}_4) = (z_1, \dot{z}_1, z_2, \dot{z}_2)$, so that the dynamics of the system is described by $n=4$ variables corresponding to the displacements and velocities of the two oscillators. It is worth stressing that, although addressing a real system is beyond the scope of this investigation, several microstructures feature nonlinear interactions between two modes, called internal resonances, that can be effectively modelled as in Eq.\eqref{eq:coupledOscillator12}, see e.g.\ the case of MEMS gyroscope analysed in \cite{art:Giorgio21}.

We have assumed to be able to observe only $y(t)=z_1(t)$. The aim is therefore to reconstruct the entire system dynamics $\mathbf{f}$ through SINDy by using only snapshots of $z_1$, thus considering realistic setting in which both $\mathbf{f}$ and $\mathbf{h}$ must be identified starting only from system observations, with the dimension of the observations smaller than the state of the system, \i.e. $o<n$. We make use of this example to show how to determine a SINDy model when the full state of the system can not be observed, and how to use this model to identify the system parameters. Such parameters can be related to unobserved system components, like the linear stiffness $k_2$ of the hidden oscillator and the coupling parameters $\alpha$ and $\beta$. Here, we have employed the time delay embedding technique, presented in Sec. \ref{sec: timedelay}, to enrich the partial observation vector and recover the original dimensionality of the system state.

We have verified that without time delay embedding, it is impossible to get a SINDy model that accurately predicts the system dynamics $\mathbf{f}$, thereby allowing system identification. In this scenario, the library $\boldsymbol{\Theta}_f$ would be solely a function of $z_1$ and the system parameters to be estimated (as we will highlight in the next subsections, we have separately addressed the estimation of $k_2$ and the coupling parameters). A SINDy model based on such a function library cannot catch some fundamental aspects of the system dynamics, such as the dependence of the damping force on $\dot{z}_1$. Considering that it is always possible to recover $\dot{z_1}$ from $z_1$, we have also attempted to train a SINDy model by using both $z_1$ and $\dot{z_1}$, but even in this case, the prediction accuracy of the model is insufficient for system identification purposes. This inadequacy persists regardless the threshold hyperparameter $L$, as it is fundamentally linked to the impossibility of including the dependence on the second oscillator. These analyses have confirmed the necessity of using time delay coordinates. After having trained a SINDy model on the new set of embedded-coordinates, we have thus proceeded with the system identification with EKF-SINDy. 

\subsubsection{Time-delay embedding}
\label{sec:timeDelayEmbedding}
In general, (fully) learning $\mathbf{f}$ directly in its original dimension from lower dimensional measurements is not possible using SINDy, as it requires knowledge of the full underlying state space variables \cite{art:Champion19}.
To tackle this issue, we adopt the time-delay embedding techniques \cite{art:Brunton17} to lift the low-dimensional time series $y(t)=z_1(t)$ into a high-dimensional space (see Sec.~\ref{sec: timedelay}).  With time-delay embedding is possible to recover hidden features carrying on information about the unobserved oscillator, which might be fundamental to describe the system's dynamics \cite{bakarji2023discovering}.


Putting the method into practice, we report the outcome of the time-delay embedding including in the Hankel matrix the dependence on the stiffness $k_2$ of the unobserved oscillator, thus setting $\phi=k_2$. These embedded coordinates will be exploited to construct the SINDy model. For the $\iota$-th sample $\phi_{\iota}$, with $\iota=1,\ldots,I$, Eqs. \eqref{eq:coupledOscillator12} are numerically integrated in time. By collecting $z_1(\phi_{\iota},t_j)$, for $j=1,\ldots,T$, 
we assemble the matrix: 
\begin{equation}
     \mathbf{A}_{\iota}=
     \begin{bmatrix}
        z_1(\phi_\iota,t_1) & z_1(\phi_\iota,t_2) & \cdots & z_1(\phi_\iota,t_{T-w}) \\
        z_1(\phi_\iota,t_2) & z_1(\phi_\iota,t_3) & \cdots & z_1(\phi_\iota,t_{T-w+1}) \\
        \vdots & \vdots & \ddots & \vdots \\
        z_1(\phi_\iota,t_w) & z_1(\phi_\iota,t_{w+1}) & \cdots & z_1(\phi_\iota,t_{T})
    \end{bmatrix}.
\label{eq:HankelMatrix_iota}%
\end{equation}

Thus, we construct the Hankel matrix $\mathbf{A}$
by stacking $\mathbf{A}=[\mathbf{A}_1,\ldots,\mathbf{A}_I]\in\mathbb{R}^{w\times (T-w)I}$,
and next we take its SVD $\mathbf{A} = \mathbf{U}\mathbf{S}\mathbf{V}^{\top},$
with $\mathbf{U}\in\mathbb{R}^{w\times w}$ and $\mathbf{V} \in\mathbb{R}^{(T-w)I\times (T-w)I}$, such that  $\mathbf{U}^\top\mathbf{U}=\mathbf{I}$ and $\mathbf{V}^\top\mathbf{V}=\mathbf{I}$, while ${\mathbf{S}}\in\mathbb{R}^{w\times (T-w)}I$ is a pseudo-diagonal matrix collecting the associated singular values $s_{\lambda}$, with $\lambda=1,\ldots,w$.

For this case, we have considered $I = 16$, $w=200$ and $\zeta=1$, where $w$ and $\zeta$ indicate the number of delay embeddings and the the number of lags in each successive embedding, respectively. The $I$ values of the parameter $k_2$ have been sampled in the range $[1.0,4.0]$ using a stratified sampling \cite{book:Saltelli07}. In practice, the range $[1.0,4.0]$ has been partitioned into $I$ subintervals, and random sampling has been performed from each subinterval. The number of parameter values required to detect parametric dependence is related to the complexity of this dependency and the problem at hand. In this specific case, being a didactic example, we chose a number that is sufficiently large for this purpose. In applications where the computation of system trajectories can be extremely costly in terms of computational resources, a more parsimonious choice of the number of instances $I$ can be considered.

\begin{figure}[b!]
  \centering
  \begin{minipage}{0.49\textwidth}
    \centering
    \includegraphics[width=\linewidth]{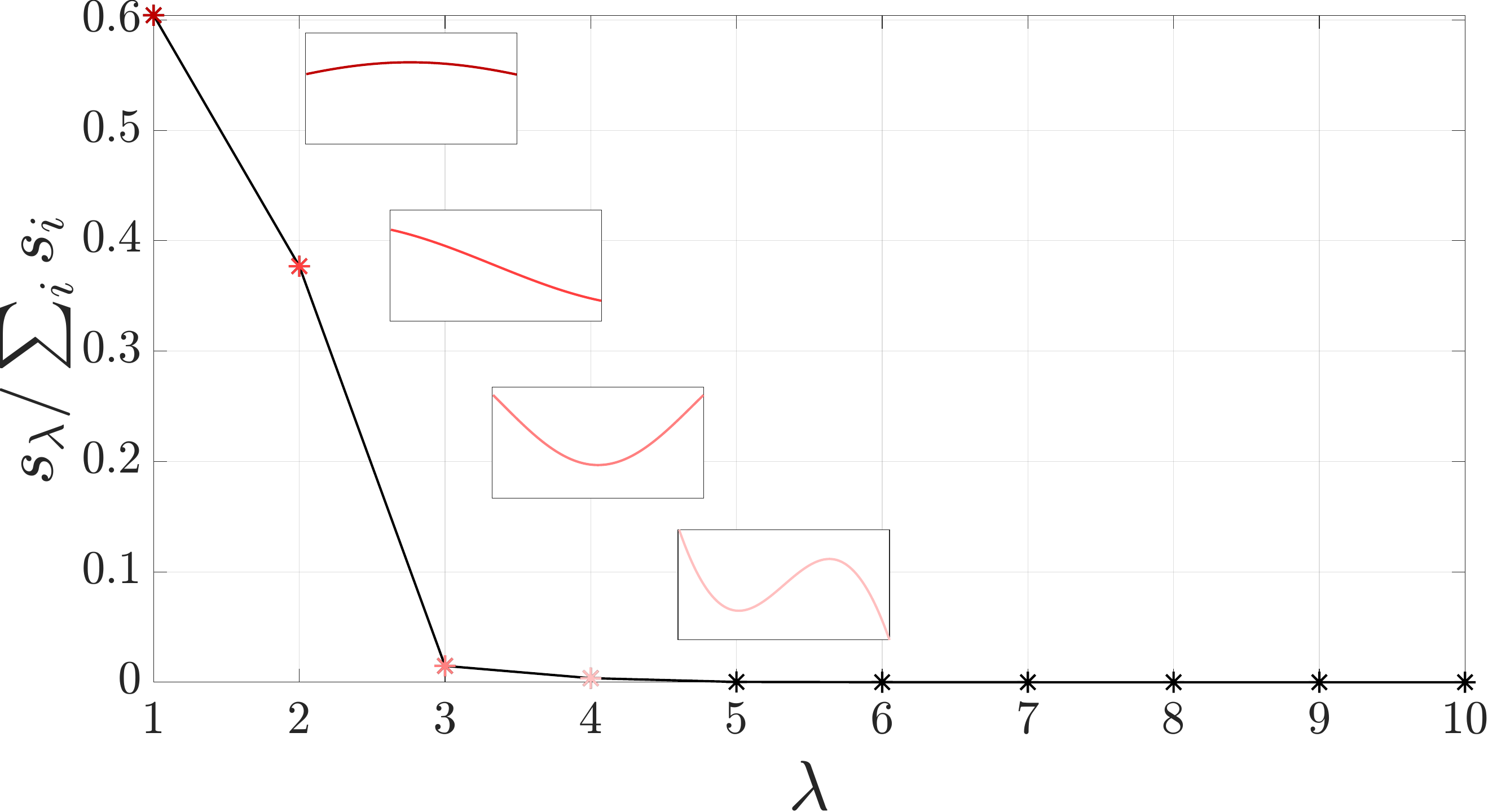}
  \end{minipage}\hfill 
  \begin{minipage}{0.49\textwidth}
    \centering
    \includegraphics[width=\linewidth]{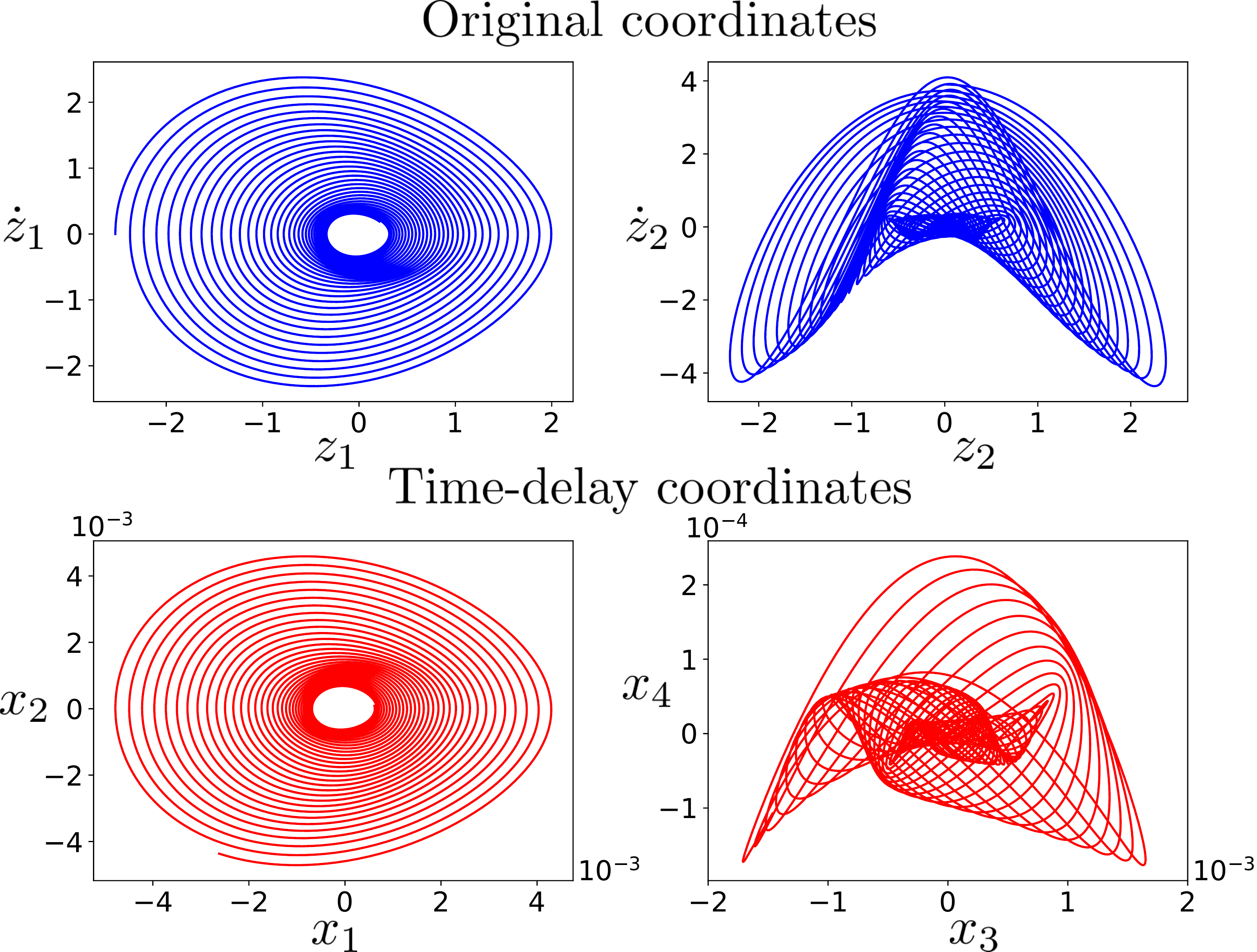} 
  \end{minipage}
{\caption{\footnotesize Nonlinear resonator. Dominant singular values $s_{\lambda}$ and left singular vectors (left) and corresponding time-delay coordinates (right), obtained by performing the SVD on the Hankel matrix constructed by integrating the considered nonlinear dynamical system for $I=16$ values of $k_2$ sampled from $[1.0,4.0]$. \\
}\label{fig:eig_1_1}}
\end{figure}

Performing the SVD on the Hankel matrix indicates that only $4$ dominant modes are responsible for the $99.97\%$ of the variance of $\mathbf{A}$, as illustrated by Fig. \ref{fig:eig_1_1}. This outcome aligns with the problem original dimension. 
In Fig. \ref{fig:eig_1_1}, we have also highlighted that the left singular values (reported in the boxes close to the corresponding singular values) resemble Legendre polynomials, as expected by choosing a small embedding period ($\zeta = 1$) \cite{vautard1989singular, broomhead1986extracting}.

Coming back to the time delay parameters $w$ and $\zeta$, we can argue that our choices guaranteed to unfold the attractor in the embedding space \cite{ma2006selection} by satisfying the Taken's embedding theorem condition on the number of delay embeddings $w>2n$. Indeed, this condition guarantees the existence of the diffeomorphism to the original system. A comparison of time-delayed coordinates with respect to the original (unobserved) ones is reported in Fig. \ref{fig:eig_1_1}. The similarity between these coordinates highlights that the time-delay coordinates adopted effectively allow for unfolding the attractor. According to the singular value decay, we have truncated the SVD matrices to account just for the $n=4$ dominant modes. The selected singular values and the left singular vectors are collected in the matrices $\tilde{\mathbf{S}}\in\mathbb{R}^{n \times n}$ and $\tilde{\mathbf{U}}\in\mathbb{R}^{w \times n}$, respectively.


\subsubsection{EKF-SINDy leveraging embedded coordinates}
\label{sec:EKF_sindy_embedding}
A SINDy model has then been trained by leveraging on the time-delayed coordinates to evolve the system dynamics. Specifically, the training has been performed after projecting the time series used to construct the Hankel matrix into the space spanned by the columns of $\tilde{\mathbf{U}}$ as in the following:
\begin{equation}
    \mathbf{x}(t_j) = \tilde{\mathbf{S}}^{-1}\tilde{\mathbf{U}}^{\top}\left( \mathbf{A} \mathbf{e}_j\right),
    \label{eq:projHankel}
\end{equation}
where $\mathbf{e}_j\in\mathbb{R}^{(T-w)I}$ is the basis vector extracting the $j$-th column of $\mathbf{A}$.

The state vector is then augmented with the corresponding parameter as $\boldsymbol{\varkappa} = \left[\mathbf{x}^{\top},\phi\right]^{\top}$.
State data are collected in the matrix $\mathbf{X}_{\varkappa}\in~\mathbb{R}^{TI\times (n+1)}$ that is necessary for training the SINDy model. According to Algorithm \ref{al:offPhase}, we use SINDy to build the function $\mathbf{f}\approx\boldsymbol{\Xi}_f\boldsymbol{\Theta}_f(\boldsymbol{\varkappa})$, describing the dynamics of the delayed coordinates. Polynomial terms up to the third order are included in $\boldsymbol{\Theta}_f$, due to the cubic nonlinearities present in the system \eqref{eq:coupledOscillator12}. Sparsity is promoted using STLSQ with $L=10^{-3}$. Interestingly, diminishing this threshold term to $L=5\cdot10^{-4}$ including more terms in the description of the dynamics, slightly improves the reconstruction capacity of SINDy but greatly deteriorates the identification capacity of the procedure. This suggests that the extra term included in the description of the system dynamics obtained for $L=5\cdot10^{-4}$ does not have a physical meaning, but simply helps in fitting the collected time histories \cite{art:Champion19}. 
On the other hand, the function $\mathbf{h}$ extracting observations $z_1(t_j)$ from the state variables $\mathbf{x}(t_j)$ is known. Specifically, $\mathbf{h}$ is linear and it can be determined by inverting Eq. \eqref{eq:projHankel} as in the following
\begin{equation}
    z_1(t_j) = \mathbf{e}^{\top}_1 \tilde{\mathbf{U}} \tilde{\mathbf{S}} \mathbf{x}(t_j),
    \label{eq:case2Jacobian}
\end{equation}
where $\mathbf{e}_1\in \mathbb{R}^w$ is the basis vector extracting the first row of $\tilde{\mathbf{U}}$ or, in other words, the first component of each left singular vector in the vector space $\mathbb{R}^w$.

Finally, to apply the EKF-SINDy based joint estimation procedure, it is necessary to compute the Jacobian matrices $\mathbf{F}$ and $\mathbf{H}$. While $\mathbf{F}$ is computed as in Eq. \eqref{eq:evJacobianSINDy}, $\mathbf{H}$ can be immediately derived by Eq. \eqref{eq:case2Jacobian}.

\subsubsection{Estimation of the stiffness $k_2$ of the hidden oscillator}

\begin{figure}[b!]
\centering
\includegraphics[width=140mm]{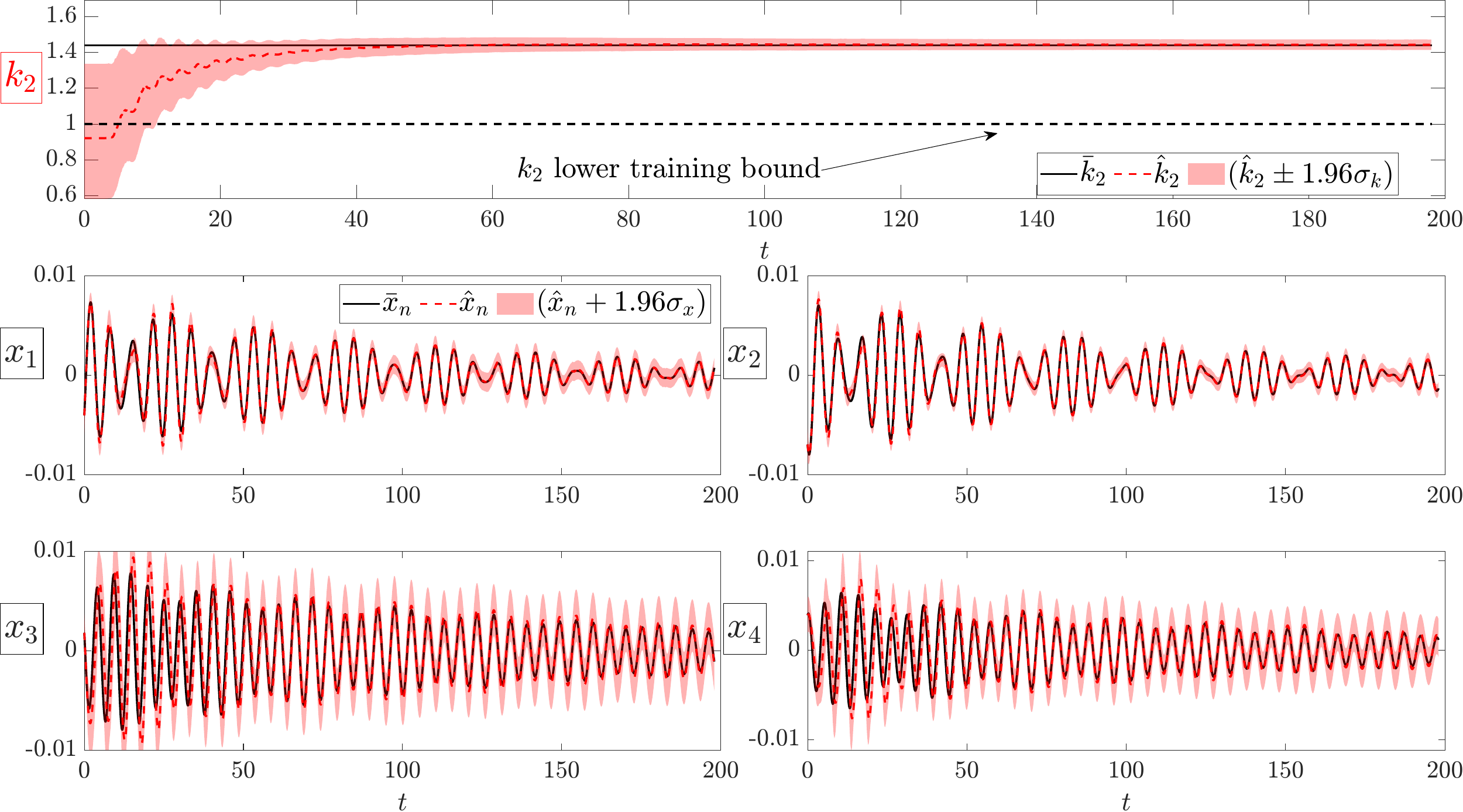}
\caption{\footnotesize 
Nonlinear dynamical system, estimation of the resonator stiffness $\bar{k}_2=1.44$. In the top graph, the estimated parameter $\hat{k}_2$ is plotted in red against the target $\bar{k}_2$ in black. The red shaded area represents the $95\%$ confidence interval of the estimates, determined using the posterior covariance. In the other graphs, the evolution of the (hidden) system response tracked by the filter is plotted in red against the real system dynamics in black.
}
\label{fig:state_1_2}
\end{figure}

We have applied the EKF-SINDy procedure to estimate the stiffness $k_2$ of the hidden oscillator. In the offline phase, we have performed time-delay embedding and approximate the dynamic model $\mathbf{f}$ by training SINDy on the time-delayed coordinates. The parameters defining the operated time delay-embedded have been reported in Sec. \ref{sec:timeDelayEmbedding}. As previously specified, $k_2$ has been sampled from $[1.0,4.0]$.

We test our method for two instances $\bar{k}_2\in\{1.44, 5.29\}$, whose results are reported in Fig.~\ref{fig:state_1_2} and Fig.~\ref{fig:state_1_4}, respectively. The tuning parameters of the filter are reported in Appendix B.

\begin{figure}[b!]
\centering
\includegraphics[width=145mm]{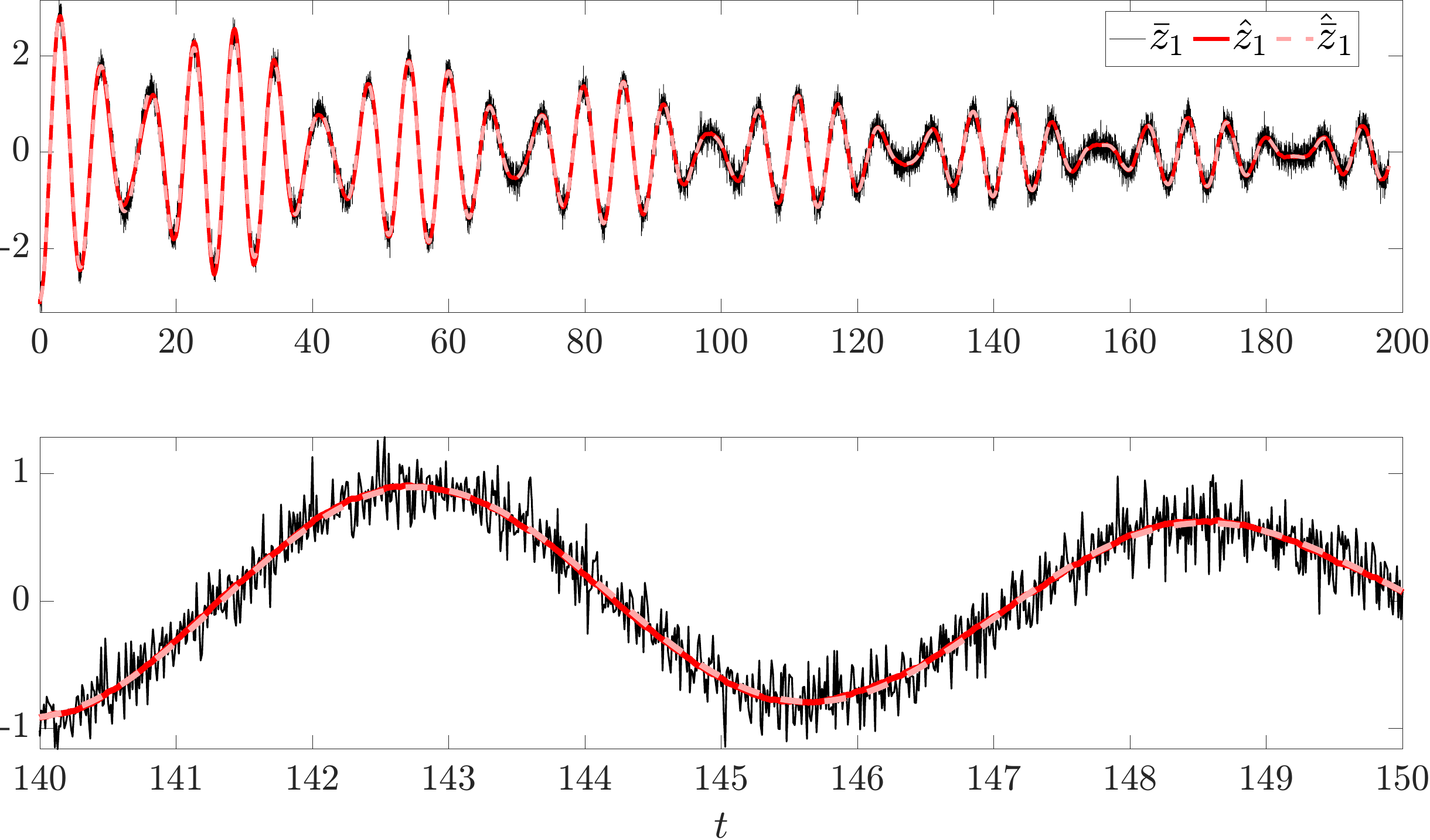}
{\caption{\footnotesize 
Nonlinear dynamical system, estimation of the resonator stiffness $\bar{k}=1.44$. Comparison between the observed variable $\bar{z}_1$ in black and its filter reconstruction $\hat{z}_1$ in red. The noise free version $\hat{\bar{z}}_1$ of the acquired measurements is depicted with a dotted line in pink.}\label{fig:out_1_2}}
\end{figure}

For the test case $\bar{k}_2 = 1.44$, starting from an initial parameter guess outside the training range of SINDy and underestimating $\bar{k}_2$ by $35\%$, the EKF-SINDy method progressively increases accuracy and decreases uncertainty until $t=50$, when it converges to the correct value. Similarly, it provides an accurate state estimate with uncertainty bounds that include the (partially observed) state of the system. We recall that data assimilation is limited to the observation of the displacement of the first oscillator, while the other quantities required to describe the system ($z_2$, $\dot{z}_1$ and $\dot{z}_2$) are unobserved. The target $\bar{\mathbf{x}} = [\bar{x}_1, \bar{x}_2, \bar{x}_3, \bar{x}_4]^\top$
is reconstructed from the real dynamics of the system for sake of comparison with the filter estimates $\hat{\mathbf{x}}_n$. After convergence at $t=~50$, the predicted mean value of the state closely matches the evolution of the state and demonstrates robustness over the entire time horizon. The reconstructed $\hat{z}_1$ is plotted against the acquired measurements $\bar{z}_1$ in Fig. \ref{fig:out_1_2}. demonstrating the noise tolerance of the procedure, with the filter estimate $\hat{z}_1$ closely aligning with the noise-free version of the signal $\hat{\bar{z}}_1$.

Starting the identification process from an initial guess for $k_2$ outside the limits for which data have been collected in the offline phase, as shown in Algorithm \ref{al:offPhase}, assesses the robustness of procedure. Operating outside such a training range is hardly achievable by other ML techniques like NNs \cite{chp:SpringerChp22,art:CAS22}, which constitutes a major advantage of the proposed procedure. This advantage stems from the capacity of the retained SINDy library terms to describe the dominant dynamical behaviour of the system. Thus, these terms approximate the physical process that underlies the observed dynamics \cite{art:Brunton16}.

To further evaluate the capability of the proposed approach to function beyond the training range of SINDy, a second test case with $\bar{k}_2=7.62$ has been conducted. In this case, the initial guess for the stiffness of the second oscillator overestimated the correct value by $20\%$. Unlike the scenario with $\bar{k}_2=1.2$, the entire identification procedure was expected to operate outside the training range $[1.0,2.0]$ of SINDy. However, EKF-SINDy successfully performed the estimation of the system properties as illustrated in Fig. \ref{fig:state_1_4}. In contrast, the estimation of the hidden state components $x_3$ and $x_4$ shows some degradation. 

\begin{figure}[h!]
\centering
\includegraphics[width=145mm]{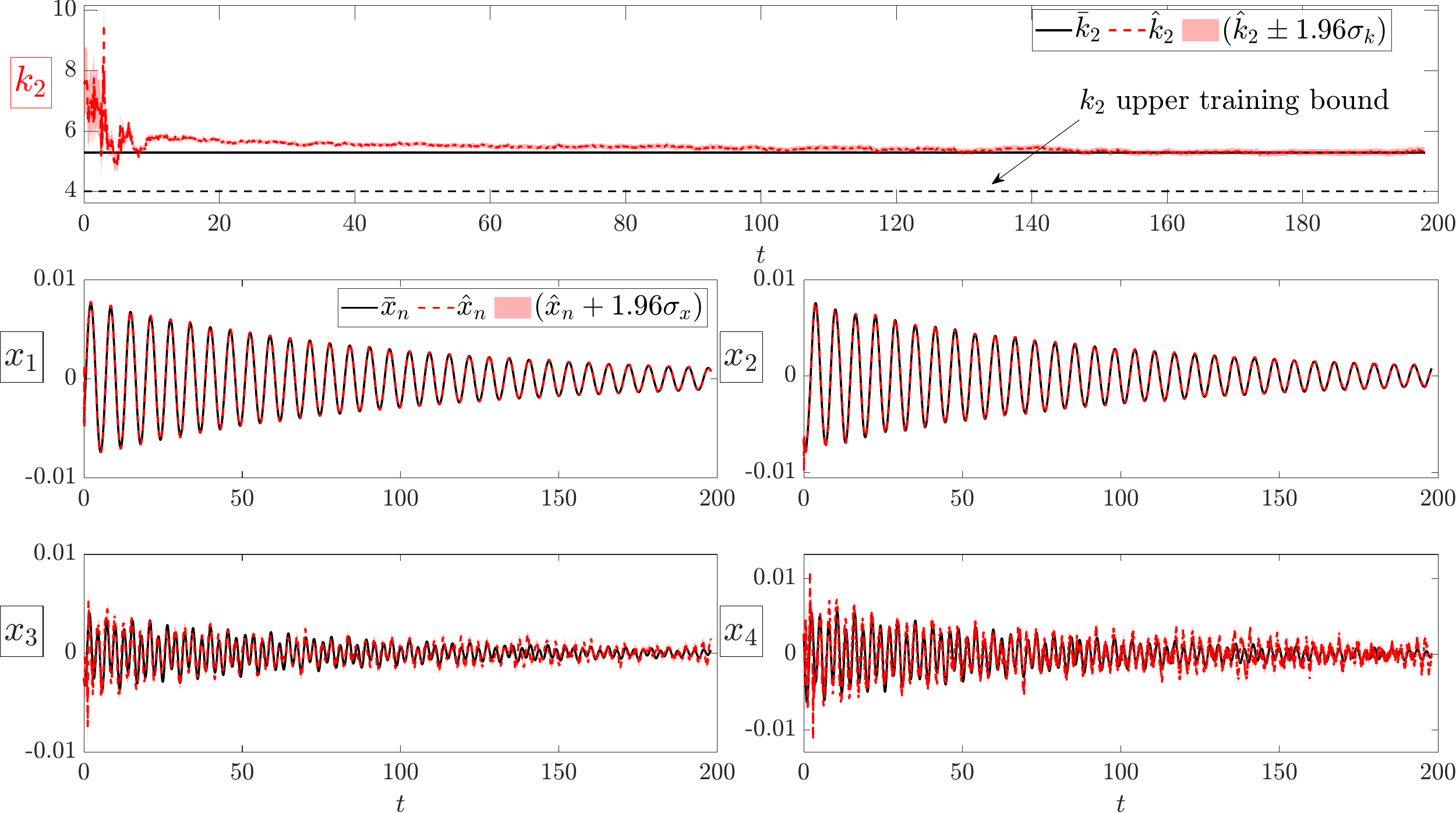}
{\caption{\footnotesize 
Nonlinear dynamical system, estimation of the resonator stiffness {\color{red}{$\bar{k}=5.29$}}. See caption of Fig. \ref{fig:state_1_2} for further details.}\label{fig:state_1_4}}
\end{figure}

\subsubsection{Estimation of the parameters $\alpha$ and $\beta$ related to linear and quadratic coupling terms}
\label{sec:couplingCoeff}

To demonstrate the procedure capability of identifying more than one system parameter, we have simultaneously addressed the estimation of the linear and quadratic coupling terms $\alpha$ and $\beta$, while setting $k_2=1.95^2 k_1$. In future work, we will investigate the internal resonances potentially arising for this ratio as done in \cite{art:Giorgio21}.

We have sampled $I=20$ values of $\alpha$ and $\beta$ to construct the Hankel matrix and to train SINDy on the time-delayed coordinates. Specifically, we have applied the Latin hypercube sampling across the domain $[0.005,0.5]^2$ defined for $\alpha$ and $\beta$ with logarithmic scaling \cite{book:Saltelli07}. The logarithmic scaling has been used to put attention to cases featuring, possibly simultaneously, low values of $\alpha$ and $\beta$. Indeed, considering only high values of $\alpha$ and $\beta$ will possibly make the system response dynamics ruled by few dominant dynamic effects, precluding the investigation of potential interaction phenomena. As in previous cases, we have fixed the number $w$ of delay embeddings to $200$, and the embedding period $\zeta$ at $1$. The tuning of the filter is reported in Appendix B.

\begin{figure}[t!]
\centering
\includegraphics[width=145mm]{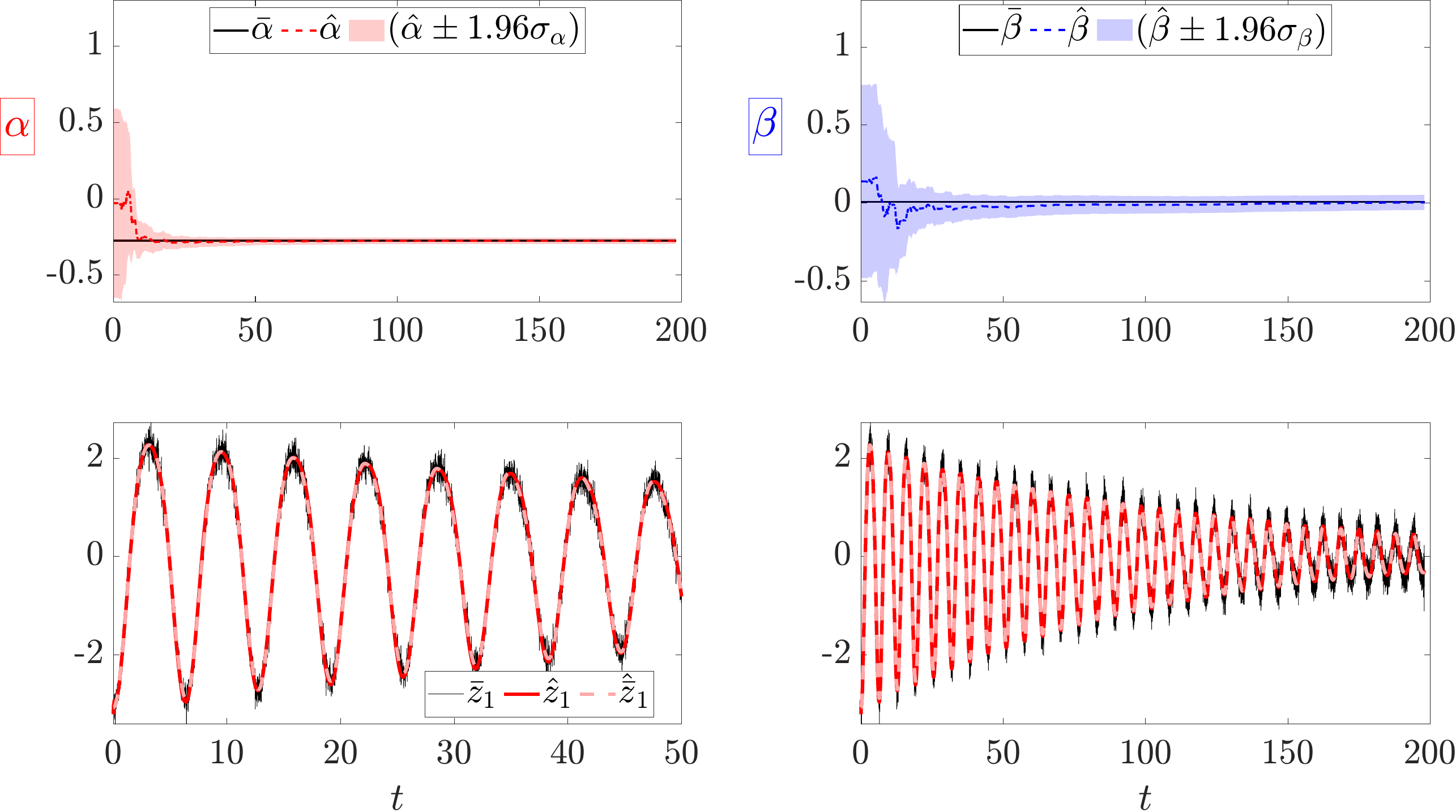}
{\caption{\footnotesize 
Nonlinear dynamical system, estimation of the parameters $\alpha$ and $\beta$ of the linear and quadratic coupling terms. In the top left corner, the estimated parameter $\hat{\alpha}$ is plotted in red against the target $\bar{\alpha}$ in black. On the right, the estimated $\hat{\beta}$ is plotted in blue against the target $\bar{\beta}$ in black. The red and blue shaded area represents $95\%$ confidence interval of the estimates, determined using the posterior covariance. In the bottom right corner, a comparison between the observed variable $\bar{z}_1$ in black and its filter reconstruction $\hat{z}_1$ in red is reported. The noise free version $\hat{\bar{z}}_1$ of the acquired measurements is depicted with a dotted line in pink. A zoom in of this plot is reported on the left.}\label{fig:2_1SENDy_param_unc}}
\end{figure}

\begin{figure}[t!]
\centering
\includegraphics[width=145mm]{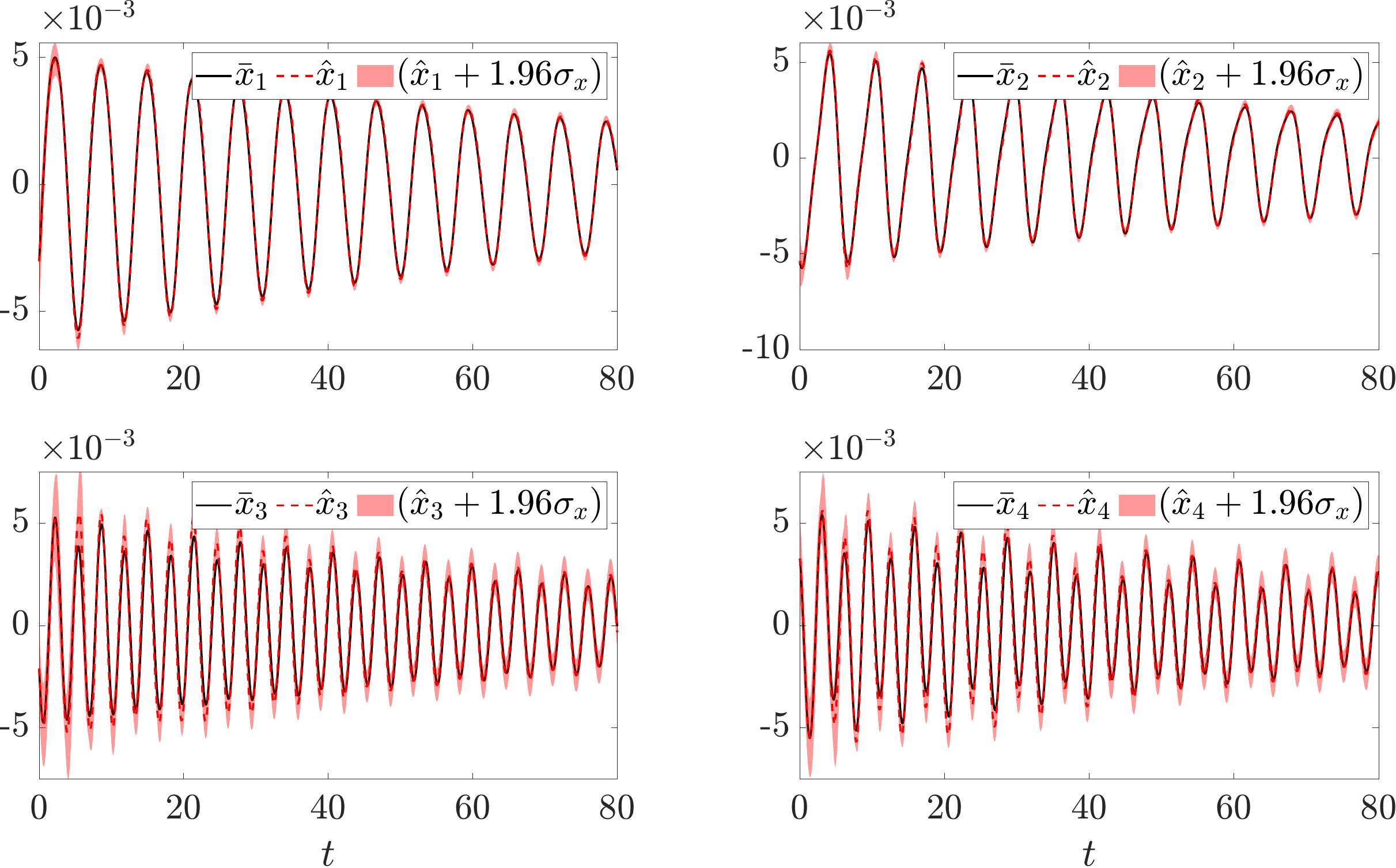}
{\caption{\footnotesize 
Nonlinear dynamical system, estimation of the parameters $\alpha$ and $\beta$ of the linear and quadratic coupling terms. The evolution of (hidden) system dynamics tracked by the filter is plotted in red against the real dynamics in black. The red shaded area represents $95\%$ confidence interval of the estimates determined using posterior covariance. For the sake of presentation, the outcomes are truncated at $t=80$.}\label{fig:2_1SENDy_state_unc}}
\end{figure}

The outcome of the system identification, illustrated in Fig. \ref{fig:2_1SENDy_param_unc}, highlights the capacity of the EKF-SINDy procedure to simultaneously estimate $\alpha$ and $\beta$. The convergence in the estimation of $\alpha$ and $\beta$ occurs roughly at $t=50$. Before that point, the prediction of the hidden dynamic state components $x_3$ and $x_4$ is imprecise, as shown in Fig. \ref{fig:2_1SENDy_state_unc}. In contrast, there is no discrepancy observed between predicted and target values for the other hidden coordinates $x_1$ and $x_2$. These observations lead to conclude that enlarging the coordinates system from the number of observed variable, namely $1$, to the number of variables, namely $4$, is crucial to for correctly describing the dynamics of the system.

\section{Conclusions}
\label{sec:conclusions}

In this work, we have proposed to empower the extended Kalman filter with the sparse identification of nonlinear dynamics (SINDy) to provide a robust, easy-to-use, and computationally inexpensive tool to enable the construction and update of digital twins. The procedure has been able to identify the dynamics of the system and their mechanical properties, providing confidence bounds for these quantities, by assimilating noisy and possibly partial observations.

We have demonstrated this capacity by addressing two case studies. First, we addressed a linear non-autonomous system consisting of a shear building model, frequently used in vibration monitoring of civil structures, excited by real seismograms. Last, we faced a partially observed nonlinear system identifying in a first analysis the stiffness of the unobserved resonator, and in a second analysis the linear and quadratic coupling coefficients of the two resonators. Interestingly, we have showed how the procedure yields good outcomes also when the identification was operated outside the training range of SINDy. This demonstrates robustness and superior generalization capabilities of the proposed method with respect to alternative machine learning strategies, which rely on fully black-box, data-driven approaches, as, e.g., neural networks. We have also shown how the time delay embedding can be used to uplift the dimensionality of the partial observations thus enabling the full description of the system dynamics.

Future and challenging directions involve applying the method to experimental data, which requires extending the approach to handle high-dimensional data. Given that SINDy is sensitive to data dimensionality, the strategy is to integrate the method with suitable dimensionality reduction techniques, such as, e.g., proper orthogonal decomposition and/or autoencoders \cite{champion2019data, bakarji2023discovering}. This integration aims to employ SINDy within a reduced space while preserving the parametric dependency as in \cite{conti2023reduced}. To further enhance the use of the procedure, we will also consider to employ automatic tuning procedures based on genetic algorithms \cite{art:Rapp03} or on swarm intelligence \cite{art:Laamari15} to perform the automatic tuning of the filter.

\section*{Code and data accessibility}
The source code of the proposed method is made available from the GitHub repository:\\
https://github.com/ContiPaolo/EKF-SINDy \cite{EKFSINDY_repo}. 

\section*{Acknowledgments}
LR and PC are supported by the Joint Research Platform ``Sensor sysTEms and Advanced Materials" (STEAM) between Politecnico di Milano and STMicroelectronics Srl. AM is member of the Gruppo Nazionale Calcolo Scientifico-Istituto Nazionale di Alta Matematica (GNCS-INdAM) and acknowledges the project “Dipartimento di Eccellenza” 2023-2027, funded by MUR, the project FAIR (Future Artificial Intelligence Research), funded by the NextGenerationEU program within the PNRR-PE-AI scheme (M4C2, Investment 1.3, Line on Artificial Intelligence), and the PRIN 2022 Project “Numerical approximation of uncertainty quantification problems for PDEs by multi-fidelity methods (UQ-FLY)” (No. 202222PACR), funded by the European Union - NextGenerationEU.
AF acknowledges the PRIN 2022 Project “DIMIN- DIgital twins of nonlinear MIcrostructures with iNnovative model-order-reduction strategies” 
(No. 2022XATLT2) funded by the European Union - NextGenerationEU.

\begin{figure}[h!]
\centerline{
\includegraphics[width=140mm]{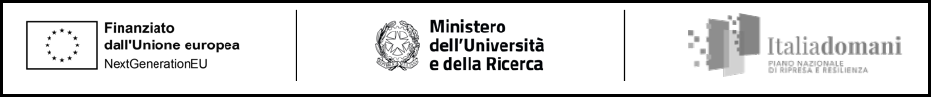}
}
\end{figure}


\section*{APPENDIX A}
\setcounter{equation}{0}
\renewcommand{\theequation}{A.\arabic{equation}}
\renewcommand{\tablename}{Table A.}

\paragraph{Sparsity promoting parameter}
The value of the threshold parameter $L$ affects the sparsity of the discovered functional relations. Conversely, it has not been necessary to adjust the regularisation strength $\delta_r$, which has been fixed at the default value of $0.05$ in the PySINDy package. In Tabs. A.1, we report the value of $L$ used for each case study.

\begin{center}
\begin{tabular}{ccc}
\toprule
\textbf{Shear building} & \textbf{Nonlinear system (}$k_2$\textbf{)} & \textbf{Nonlinear system} \text{ (}$\alpha$ \text{ and }$\beta$\text{)} \\
\midrule
$10^2$ & $5\cdot 10^{-4}$ & $10^{-3}$ \\
\bottomrule
\end{tabular}
\captionof{table}{\footnotesize Sparsity promoting threshold parameter $L$.}
\end{center}

The resulting functional relations are detailed in the following, specifying the order of the polynomial terms (linear, quadratic or cubic) employed by the trained SINDy models, and the state variables relevant to each case study.

\paragraph{Shear building under seismic excitations}
The state vector $\boldsymbol{\varkappa}$ collects the floor displacements, velocities, and interstorey stiffness $k$, namely $\boldsymbol{\varkappa}=[x_1,x_2,\dot{x}_1,\dot{x}_2,k]^{\top}$. In this way, SINDy can be directly applied to model the evolution of the state vector $\boldsymbol{\varkappa}$ as a first order differential equation, where state variables have been suitably rescaled to have comparable magnitude. Polynomial terms up to the second order are considered for the state variables, while extra linear terms are used to account for the external forcing that makes the system non-autonomous. Specifically, the forcing terms $b_1$ and $b_2$ are applied respectively to the first and second interstoreys.  

The identified SINDy model is reported here at the second order consistently with the form of Eq. \eqref{eq:shearBuilding} ruling the dynamics of the system. For the first order formulation, the reader can refer to the GitHub repository \cite{EKFSINDY_repo}.

\begin{subequations}
\begin{align}
\begin{split}
\qquad \qquad \qquad \qquad \qquad \qquad \qquad \qquad \qquad \qquad \qquad \qquad \qquad \qquad 
\ddot{x}_1 = \\
-0.716 x_1 + 0.446 x_2 -2.25 \dot{x}_1 + 1.112 \dot{x}_2 & \qquad \text{(linear)}\\
+ 308259736.25 b_1 -308259740.5 b_2 & \qquad \text{(forcing)}\\
-0.088 {x_1}^2 + 0.055 x_1 x_2 + 2.984 x_1 \dot{x}_1 -1.883 x_1 \dot{x}_2 -351.04 x_1 k  -1.844 x_2 \dot{x}_1 & \qquad \text{(quadratic)}\\
+ 1.164 x_2 \dot{x}_2 + 174.107 x_2 k -0.534 {\dot{x}_1}^2 + 0.598 \dot{x}_1\dot{x}_2 -3.768 \dot{x}_1 k -0.164 {\dot{x}_2}^2 + 1.776 \dot{x}_2 k & \qquad \text{(quadratic)}
\end{split}\\
\begin{split}
\ddot{x}_2 = \\
-0.557 x_1 + 0.346 x_2 + 0.922 \dot{x}_1 -1.023 \dot{x}_2 & \qquad \text{(linear)}\\
+ 307191775.013 b_1 -307191779.264 b_2 & \qquad \text{(forcing)}\\
-3.33 {x_1}^2 + 4.031 x_1 x_2 -0.01 x_1 \dot{x}_1 + 193.088 x_1 k -1.218 {x_2}^2 + 0.056 x_2 \dot{x}_1 & \qquad \text{(quadratic)}\\
-0.03 x_2 \dot{x}_2 -188.664 x_2 k -0.09 {\dot{x}_1}^2 + 0.036 \dot{x}_1 \dot{x}_2 + 3.186 \dot{x}_1 k + 0.014 {\dot{x}_2}^2 -2.861 \dot{x}_2 k & \qquad \text{(quadratic)}
\end{split}
\end{align}
\end{subequations}

\paragraph{Nonlinear dynamical system, estimation of the stiffness $k_2$ of the hidden oscillator}

The state vector $\boldsymbol{\varkappa}$ collects the time delay coordinates and the $k_2$, namely $\boldsymbol{\varkappa}=[x_1,x_2,x_3,x_4,k_2]^{\top}$. Considering a polynomial terms up to the third order, the total number of library functions is $55$. According to the set threshold level $L$, SINDy employs: $4$ terms to model $\dot{x}_1$; $5$ terms for $\dot{x}_2$; $7$ terms for $\dot{x}_3$; $9$ terms for $\dot{x}_4$.

\begin{subequations}
\begin{align}
\begin{split}
\qquad \qquad \qquad \qquad \qquad \qquad \qquad \qquad \qquad \qquad \qquad \qquad \qquad \qquad \qquad \qquad \dot{x}_1 = \\
-0.007 x_1 -1.027 x_2 -0.001 x_3 -0.009 x_4 & \qquad \text{(linear)}
\end{split}\\
\begin{split}
\dot{x}_2 = \\
1.004 x_1 -0.020 x_2 -0.150 x_3 & \qquad \text{(linear)}\\
-0.001 x_1 k_2 + 0.004 x_3 k_2 & \qquad \text{(quadratic)}
\end{split}\\
\begin{split}
\dot{x}_3 = \\
0.005 x_1 + 0.157 x_2 -0.002 x_3 -1.246 x_4 & \qquad \text{(linear)}\\
-0.001 x_1 k_2  -0.006 x_2 k_2 + 0.019 x_4 k_2 & \qquad \text{(quadratic)}
\end{split}\\
\begin{split}
\dot{x}_4 = \\
0.262 x_1 -0.049 x_2 -0.061 x_3 -0.011 x_4 & \qquad \text{(linear)}\\
-0.189 x_1 k_2 + 0.044 x_2 k_2 + 0.834 x_3 k_2 & \qquad \text{(quadratic)}\\
-0.001 x_2 k_2^2 + 0.002 x_3 k_2^2 & \qquad \text{(cubic)}
\end{split}
\end{align}
\end{subequations}

\paragraph{Nonlinear dynamical system, estimation of the coupling parameters $\alpha$ and $\beta$}

The state vector $\boldsymbol{\varkappa}$ collects the time delay coordinates and the coupling parameters, namely $\boldsymbol{\varkappa}=[x_1,x_2,x_3,x_4,\alpha,\beta]^{\top}$. Considering a polynomial terms up to the third order, the total number of library functions is $83$. According to the set threshold level $L$, SINDy employs: $9$ terms to model $\dot{x}_1$; $12$ terms for $\dot{x}_2$; $14$ terms for $\dot{x}_3$; $29$ terms for $\dot{x}_4$.

\begin{subequations}
\begin{align}
\begin{split}
\dot{x}_1 = \\
-0.006 x_1 -1.013 x_2 -0.001 x_3 -0.034 x_4 & \qquad \text{(linear)} \\
+ 0.001 x_4\beta & \qquad \text{(quadratic)}\\
-0.001 x_2 \alpha^2 -0.001 x_2 \alpha \beta + 0.001 x_4 \alpha^2 + 0.001 x_4 \alpha \beta & \qquad \text{(cubic)}
\end{split} \\
\begin{split}
\dot{x}_2 = \\
+0.996 x_1 -0.014 x_2 -0.343 x_3 & \qquad \text{(linear)}\\
+ 0.002 x_1 \alpha + 0.002 x_3 \alpha + 0.005 x_3 \beta & \qquad \text{(quadratic)}\\
-0.028 x_1 \alpha^2 + -0.007 x_1 \alpha \beta -0.004 x_1 \beta^2 + 0.007 x_3 \alpha^2 + 0.021 x_3 \alpha \beta + 0.011 x_3 \beta^2 & \qquad \text{(cubic)}
\end{split} \\
\begin{split}
\dot{x}_3 = \\
-0.003 x_1 + 0.338 x_2 -0.001 x_3 -1.924 x_4 & \qquad \text{(linear)} \\
+ 0.002 x_1 \alpha -0.002 x_2 \beta + 0.011 x_4 \alpha + 0.015 x_4 \beta & \qquad \text{(quadratic)}\\
-0.042 x_2 \alpha^2 -0.021 x_2 \alpha \beta -0.012 x_2 \beta^2 + 0.029 x_4 \alpha^2 + 0.059 x_4 \alpha \beta + 0.030 x_4 \beta^2 & \qquad \text{(cubic)}
\end{split}\\
\begin{split}
\dot{x}_4  = \\
+0.114 x_1 -0.004 x_2 + 1.849 x_3 + 0.029 x_4 & \qquad \text{(linear)}\\
+ 8.899 x_1^2 -0.422 x_1 x_4 + 0.089 x_1 \alpha + 0.035 x_1 \beta + 0.290 x_2 x_3 -0.008 x_2 \alpha & \qquad \text{(quadratic)}\\ -0.478 x_3 \alpha -0.285 x_3 \beta + 0.916 x_4^2 + 0.497 x_4 \alpha -0.003 \alpha^2 -0.001 \beta^2 & \qquad \text{(quadratic)}\\
-1.242 x_1 \alpha^2 -0.244 x_1 \alpha \beta -0.230 x_1 \beta^2 + 0.086 x_2 \alpha^2 -0.202 x_2 \alpha \beta & \qquad \text{(cubic)} \\
-0.759 x_3 \alpha^2 -0.460 x_3 \alpha \beta + 1.501 x_3 \beta^2 + 0.950 x_4 \alpha^2 & \qquad \text{(cubic)} \\
-0.392 x_4 \alpha \beta -0.006 \alpha^3 + 0.018 \alpha^2 \beta + 0.002 \beta^3 & \qquad \text{(cubic)}
\end{split}
\end{align}
\end{subequations}

\section*{APPENDIX B}
\label{sec:appendixTuning}

\renewcommand{\tablename}{Table B.}
We report further information about the tuning of the EKF for the considered cases of study. In all cases, we have initialised the state covariance matrix $\mathbf{P}(0)$, the process noise covariance matrix $\mathbf{Q}$, and the observation noise covariance matrix $\mathbf{R}$ as diagonal matrices. In Tabs. B.1-B.3, the diagonal terms of these matrices are gathered.\\

\scriptsize{
\begin{minipage}[h]{0.28\textwidth}
\label{tab:EKFtuning_1}
\begin{tabular}{lr}
\toprule
    \multicolumn{2}{c}{\textbf{Shear building}} \\
    \midrule
    $x_1, x_2, \dot{x}_1, \dot{x}_2$ & $10^{-5}$ \\
    $k$ & $2\cdot10^{-5}$\\
    \midrule
    \multicolumn{2}{c}{\textbf{Nonlinear dynamic}}\\
    \multicolumn{2}{c}{\textbf{system}}\\
    \midrule
    \multicolumn{2}{c}{$\bar{k}_2=1.2$}\\
    $x_1, x_2, x_3, x_4$ & $10^{-6}$\\
    $k_2$ & $10^{-2}$\\
    \midrule
    \multicolumn{2}{c}{$\bar{k}_2=2.3$}\\
    $x_1$ & $10^{-6}$\\
    $x_2$ & $10^{-6}$\\
    $x_3$ & $10^{-6}$\\
    $x_4$ & $10^{-6}$\\
    $k_2$ & $2\cdot 10^{-6}$\\
    \midrule
    \multicolumn{2}{c}{$\alpha$, $\beta$}\\
    $x_1$ & $10^{-6}$\\
    $x_2$ & $10^{-6}$\\
    $x_3$ & $10^{-6}$\\
    $x_4$ & $10^{-6}$\\
    $\alpha$ & $10^{-1}$\\
    $\beta$  & $10^{-1}$\\
    \bottomrule
\end{tabular}
\captionof{table}{Initial covariance $\mathbf{P}(0)$, diagonal terms. In the left column, we specify the augmented state quantities to which these values are related.}
\end{minipage}%
$\quad \quad$
\begin{minipage}[h]{0.28\textwidth}
\label{tab:EKFtuning_2}
\begin{tabular}{lr}
\toprule
    \multicolumn{2}{c}{\textbf{Shear building}} \\
    \midrule
    $x_1, x_2, \dot{x}_1, \dot{x}_2$ & $10^{-4}$\\
    $k$ & $10^{-8}$\\
    \midrule
    \multicolumn{2}{c}{\textbf{Nonlinear dynamic}}\\
    \multicolumn{2}{c}{\textbf{system}}\\
    \midrule
    \multicolumn{2}{c}{$\bar{k}_2=1.2$}\\
    $x_1, x_2, x_3, x_4$ & $10^{-7}$\\
    $k_2$ & $10^{-7}$\\
    \midrule
    \multicolumn{2}{c}{$\bar{k}_2=2.3$}\\
    $x_1$ & $10^{-10}$\\
    $x_2$ & $10^{-10}$ \\
    $x_3$ & $10^{-8}$ \\
    $x_4$ & $10^{-8}$\\
    $k_2$ & $2\cdot 10^{-6}$\\
    \midrule
    \multicolumn{2}{c}{$\alpha$, $\beta$}\\
    $x_1$ & $5\cdot 10^{-10}$\\
    $x_2$ & $5\cdot 10^{-10}$ \\
    $x_3$ & $8\cdot 10^{-8}$ \\
    $x_4$ & $8\cdot 10^{-8}$\\
    $\alpha$ & $10^{-11}$\\
    $\beta$ & $10^{-11}$\\
    \bottomrule
\end{tabular}
\captionof{table}{Process noise $\mathbf{Q}$, diagonal terms. In the left column, we specify the augmented state quantities to which these values are related.}
\end{minipage}%
$\quad \quad$
\noindent
\begin{minipage}[h]{0.28\textwidth}
\label{tab:EKFtuning_3}
\begin{tabular}{cr}
\toprule
    \multicolumn{2}{c}{\textbf{Shear building}} \\
    \midrule
    $x_1$ & $5 \cdot 10^{-1}$ \\
    $x_2$ & $5 \cdot 10^{-1}$ \\
    $\dot{x}_1$ & $5 \cdot 10^{-3}$ \\
    $\dot{x}_2$ & $5 \cdot 10^{-3}$ \\
    $\ddot{x}_1$ & $5$ \\
    $\ddot{x}_2$ & $5$ \\
    \midrule
    \multicolumn{2}{c}{\textbf{Nonlinear dynamic}}\\
     \multicolumn{2}{c}{\textbf{system}}\\
    \midrule
    \multicolumn{2}{c}{$\bar{k}_2=1.2$}\\
    $z_1$ & $10^{-7}$\\
    \midrule
    \multicolumn{2}{c}{$\bar{k}_2=2.3$}\\
    $z_1$ & $5\cdot10^{-4}$\\
    \midrule
    \multicolumn{2}{c}{$\alpha$, $\beta$}\\
    $z_1$ & $3\cdot10^{-1}$\\
    \bottomrule
\end{tabular}
\captionof{table}{Measurement noise $\mathbf{R}$, diagonal terms. In the left column, we specify the observed quantities to which these values are related.}
\end{minipage}%
}

\normalsize \bigskip

As mentioned, Kalman filter tuning has been operated through a trial-and-error procedure. Some practical guidelines can be used to expedite this procedure. These guidelines are derived by experience, and account for the probabilistic meaning of the process $\mathbf{Q}$ and observation $\mathbf{R}$ noise covariance matrices. The following heuristic discussion must thus be intended more as a practical guidance for the reader, than as a systematic approach to Kalman filter tuning.

First, it is advisable to set the entries of $\mathbf{R}$ based on the specifications of the acquisition system. Next, the identification outcome should be checked for a random initialisation of $\mathbf{Q}$. If the filter is unable to update the mean of the model parameter predictions $\hat{\boldsymbol{\phi}}$ during the acquisition window, the entries of $\mathbf{R}$ are probably too high. If the value of $\hat{\boldsymbol{\phi}}$ is instead sensitive to the acquisitions, attention should then be directed to $\mathbf{Q}$. Examining again the identification outcomes, if the confidence intervals do not reduce during the analysis, it is likely that the entries of $\mathbf{Q}$ are set too high, causing the filter to overly trust measurements and neglect model predictions. This precludes the filter from gaining confidence in its predictions, and thus reducing the confidence intervals. Therefore, it is advisable to lower the process noise parameters. In doing that, one should initially operate large changes in the entries of $\mathbf{Q}$ (even by an order of magnitude), and perform a fine tuning in the later trials. However, excessive reduction of the process noise parameters may cause divergence in the filter estimates.

\scriptsize \bigskip

\bibliographystyle{elsarticle-num}
\bibliography{references}

\end{document}